\documentclass{article}[12pt]

\usepackage{epsfig}
\usepackage{graphicx}
\usepackage{color}
\usepackage{dcolumn}
\usepackage{bm}
\usepackage{amsmath}
\usepackage{amssymb}
\usepackage{amsxtra}
\usepackage{cancel}
\usepackage{multirow}
\usepackage{undertilde}
\usepackage{calc}

\newtheorem{theorem}{Theorem}[section]
\newtheorem{definition}[theorem]{Definition}
\newtheorem{cor}[theorem]{Corollary}
\newtheorem{lem}[theorem]{Lemma}
\newtheorem{prop}[theorem]{Proposition}
\newtheorem{rem}[theorem]{Remark}

\renewcommand{\setminus}{-}

\def\bnprf{\noindent {\bf Proof} \ }
\def\edprf{$_{\Box}$}
\def\ubar#1{\underbar{$#1$}}
\def\Bar{\overline}

\def\C{\mathcal{C}^{\pm}}
\def\Cp{\mathcal{C}^{+}}
\def\Cm{\mathcal{C}^{-}}
\def\U{\mathcal{U}^{\pm}}
\def\Z{\mathcal{Z}^{\pm}}
\def\I{\mathcal{C}^{0}}
\def\gS{\mathcal{S}}
\def\gR{\mathcal{G}}
\def\gT{\mathcal{P}}
\def\gD{\mathcal{D}}
\def\gQ{\mathcal{Q}}
\def\dummyM{M}
\def\sM{\mathcal{\dummyM}}
\def\Mp{\dummyM^{+}}
\def\Mm{\dummyM^{-}}
\def\hMp{\hat{\dummyM}^{+}}
\def\hMm{\hat{\dummyM}^{-}}
\def\sK{\mathcal{K}}
\def\Kp{K^{+}}
\def\Km{K^{-}}
\def\HH{\Delta}
\def\mod{\mathrm{mod}}
\def\CC{\vec{C}_{3}}
\def\L{\mathcal{I}}
\def\Linv{\L^{-1}}
\def\IV#1{\Linv(#1)}
\def\VI#1{\L(#1)}
\def\de{arc}	
\def\DG{\vec{G}}
\def\DH{\vec{H}}
\def\DS{S}
\def\DE{A}
\def\DP{P}
\def\coloneqq{\mathrel{\mathop:}=}
\def\ngeq{\succeq}
\def\ngrt{\succ}
\def\nleq{\preceq}

\def\pgeq{\geq}
\def\pgrt{>}

\newsavebox{\bsq}
\savebox{\bsq}(6,6){{\scriptsize $\bullet$}}
\newsavebox{\sq}
\savebox{\sq}(6,6){{\scriptsize $\circ$}}

\title{Algorithms for realizing degree sequences of directed graphs}

\author{M.~Drew LaMar\thanks{Department of Applied Science,
The College of William and Mary,
311 McGlothlin-Street Hall,
Williamsburg VA 23187 ({\tt mdlama@wm.edu}).}}

\begin{document}

\maketitle

\begin{abstract}
The Havel-Hakimi algorithm for constructing realizations of degree sequences for undirected graphs has been used extensively in the literature.  A result by Kleitman and Wang extends the Havel-Hakimi algorithm to degree sequences for directed graphs.  In this paper we go a step further and describe a modification of Kleitman and Wang's algorithm that is a more natural extension of Havel-Hakimi's algorithm, in the sense that our extension can be made equivalent to Havel-Hakimi's algorithm when the degree sequence has equal in and out degrees and an even degree sum.  We identify special degree sequences, called directed 3-cycle anchored, that are ill-defined for the algorithm and force a particular local structure on all directed graph realizations.  We give structural characterizations of these realizations, as well as characterizations of the ill-defined degree sequences, leading to a well-defined algorithm.
\end{abstract}




\section{Introduction}

All graphs (digraphs) in this article will be simple, i.e. with no self-loops or multi-edges (-arcs).  Given an integer sequence $d$, we say $d$ is {\it graphic} if there exists an undirected graph $G$ with degree sequence $d$.  We say $G$ {\it realizes} $d$ and denote the set of all realizations of $d$ by $R(d)$.  Can we give any structural information on the undirected graphs in $R(d)$?  Two fundamental questions include determining when $R(d)$ is nonempty, and if so how to construct a graph in $R(d)$.  Answers to the first question include checking the Erd\H{o}s-Gallai inequalities \cite{Erdos:1961p5414}, while the Havel-Hakimi algorithm \cite{Havel:1955p7653, hakimi:496, hakimi:135} answers both the first and the second question.

In the case of directed graphs, we consider integer-pair sequences $d=\{(d_{i}^{+},d_{i}^{-})\}_{i=1}^{N}$.  Similar to above, we say $d$ is {\it digraphic} if there exists a digraph (i.e. directed graph) with degree sequence $d$, also denoting the set of digraph realizations of $d$ by $R(d)$.  When $d$ is digraphic, then $d^{+}$ and $d^{-}$ denote the out-degree and in-degree sequences of $d$, respectively.  We have similar results for directed graphs regarding the two questions above, in particular the Kleitman and Wang algorithm \cite{Kleitman:1973p4546} for constructing directed graphs from integer-pair sequences and the Fulkerson inequalities \cite{Fulkerson:1960p3387} to check existence.  It is instructive to compare and contrast the theorems and techniques from both case studies, since undirected graphs are specific examples of directed graphs when undirected edges are identified with bidirectional arcs.  

How related are the techniques, and what can their similarities and differences tell us?  This paper explores this question in regards to constructing graphs from degree sequences.  In particular, its main study is extensions of the Havel-Hakimi algorithm to directed graphs similar in nature to Kleitman and Wang's algorithm.  There is a difference between these two algorithms that is addressed, leading to a modified algorithm that is in a way more similar to Havel-Hakimi on undirected graphs.  The similarity deals with the natural identification of undirected edges with bidirectional arcs.  

The new algorithm, however, is not immediately applicable, since it can be quickly shown to be ill-defined.  In other words, there are sequences where the new algorithm fails to produce a digraph realization even if one exists.  However, we show these ill-defined sequences are highly structured, having a degree sequence characterization so that they can be identified.  Once they are identified, however, we need to know how to make the appropriate connections, so we show as well that their digraph realizations have a structural characterization.  Combining this all together leads to a well-defined algorithm.  While this algorithm is not necessarily more efficient than just performing Kleitman and Wang's algorithm, the theorems and techniques from its proof are shown to be helpful in answering a question regarding Eulerian sequences (i.e. sequences with equal in and out degrees), with the possibility that it will shed some light on other questions as well.

One particular example of a useful result from this exploration is that the ill-defined sequences are in fact shown to be equivalent to what we are calling directed 3-cycle anchored sequences, or {\it $\CC$-anchored}, where $\CC$ denotes a directed 3-cycle with vertex set $\{v_{1},v_{2},v_{3}\}$ and arc set $\{(v_{1},v_{2}),(v_{2},v_{3}),(v_{3},v_{1})\}$.  $\CC$-anchored sequences force a particular {\it local} structure on all digraph realizations.  In particular, $\CC$-anchored sequences are forcibly $\CC$-digraphic with the added local constraint that there is a set $\mathcal{J}$ of indices for the degree sequence such that for every index $i \in \mathcal{J}$ and every digraph realization of the degree sequence, there is an induced subgraph isomorphic to $\CC$ containing the labeled vertex $v_{i}$.


Section \ref{sec:digraphic} reviews some of the existing theory and tools regarding graphic and digraphic sequences.  The modification of Kleitman and Wang's algorithm is introduced here, along with a definition of the ill-defined sequences and their digraph realizations, with some results on what can be quickly deduced from these definitions.  The beginning of Section~\ref{sec:characterizations} summarizes the rest of the paper with an outline of the proof of the main results.  It is in this section where we prove that the ill-defined sequences and their digraph realizations have a degree sequence and structural characterization, respectively.  Section \ref{sec:anchored} introduces a general definition of $\DH$-anchored sequences and proves the equivalence of the ill-defined sequences with $\CC$-anchored sequences.  Finally, Section \ref{sec:algorithm} revisits the context of the algorithm itself and makes more explicit the steps necessary to guarantee a well-defined algorithm.  We end by using the algorithm to construct realizations of Eulerian degree sequences with maximal bidirectional arcs.

\subsection{Notation} \label{sec:notation}

We denote undirected graphs as $G$, where $V(G)$ is the vertex set and $E(G)$ the edge set.  Directed graphs are similarly denoted by $\DG$, with $V(\DG)$ the vertex set and $\DE(\DG)$ the {\de } set.  The remaining notation will be defined in reference to either graphs or digraphs, with a similar definition applying to the other.  The set of directed paths in digraphs will be denoted by $\DP(\DG)$.  We will drop the reference to $\DG$ when the digraph is understood through the notation $\DG = (V,\DE)$, for example.

An edge or {\de } between vertices $a$ and $b$ will be denoted by $(a,b)$, with the orientation given by the ordering in the case of directed graphs.  If $X,Y \subseteq V$, then $(X,Y) \subseteq \DE$ is defined as the set of {\de}s $(x,y) \in \DE$ such that $x \in X$ and $y \in Y$.  Given vertex sets $X,Y,Z \subseteq V$, we define in a similar manner the set of directed 3-paths $(X,Y,Z) \subseteq \DP$ as the directed 3-paths $(x,y,z) \in \DP$ such that $x \in X$, $y \in Y$ and $z \in Z$.

Given a digraph $\DG = (V,\DE)$ and vertex sets $X,Y\subseteq V$, we define the subgraph $\DG[X,Y] = (X\cup Y, \DE[X,Y])$, where $\DE[X,Y] = \{(x,y)\in \DE : x \in X \ \mbox{and} \ y \in Y\}$.  When $X=Y$, we have the usual definition of an induced subgraph and will denote this by $\DG[X]$.

In figures and diagrams, we will frequently use arrows to denote relations between vertices in the following manner:
\begin{eqnarray*}
x \rightarrow y     & \Longleftrightarrow & (x,y) \in \DE \\
x \leftarrow y      & \Longleftrightarrow & (y,x) \in \DE \\
x \leftrightarrow y & \Longleftrightarrow & \{(x,y), (y,x)\} \subseteq \DE \\
x \cdots y          & \Longleftrightarrow & \{(x,y), (y,x)\} \subseteq \DE^{C}.
\end{eqnarray*}
We will also use a dashed-dotted directional line between $x$ and $y$ to denote an allowable {\de }, i.e. $(x,y) \in \DE$ or $(x,y) \notin \DE$.  When vertices are substituted by sets of vertices, we will use the same convention above, e.g. $X \cdots Y$ if and only if $\{(x,y), (y,x)\} \subseteq \DE^{C}$ for all $x\in X$ and $y\in Y$.

All integer and integer-pair sequences will be non-negative, unless explicitly stated otherwise.  Since we will be moving back and forth between vertices, degree sequences and indices, we define a bijective index function $\L: V \longrightarrow \{1,\ldots,|V|\}$ going from vertices to indices of the degree sequence.
We will often refer to an integer-pair sequence $d=\{(d_{i}^{+},d_{i}^{-})\}_{i=1}^{N}$ in matrix form as
\[d = \left(
\begin{array}{ccc}
d_{1}^{+} & \cdots & d_{N}^{+} \\
d_{1}^{-} & \cdots & d_{N}^{-}
\end{array}
\right)
\]
with the first and second rows corresponding to $d^{+}$ and $d^{-}$, respectively.

\section{Realizing graphic and digraphic sequences} \label{sec:digraphic}

The following result by Erd\H{o}s and Gallai addresses when integer sequences are graphic:
\begin{theorem}[Erd\H{o}s-Gallai \cite{Erdos:1961p5414}]
\label{thm:erdos}
Let $d$ be a non-increasing integer sequence.  Then $d$ is graphic if and only if $\sum_{i=1}^{N} d_{i}$ is even and for $k=1,\ldots,N$,
\[k(k-1) + \sum_{i=k+1}^{N} \min\{d_{i},k\} \geq \sum_{i=1}^{k} d_{i}.\]
\end{theorem}
Note that not all inequalities are necessary, as Tripathi and Vijay \cite{Tripathi:2003p3520} have shown that you only need to check as many inequalities as there are distinct terms in the sequence.

While this establishes existence of a realization of an integer sequence, it does not give an algorithm for computing one.  The following constructive algorithm was given independently by Havel and Hakimi \cite{Havel:1955p7653, hakimi:496, hakimi:135}.  Let $d$ be an integer sequence of length $N \geq 2$.  Choose an index $i$ of $d$ such that $d_{i}\neq 0$.  If there are not $d_{i}$ positive entries in $d$ other than at $i$, then $d$ is not graphic, so suppose there are $d_{i}$ positive entries in $d$ other than at $i$.  Construct the {\it residual degree sequence} $\hat{d}$ by setting $d_{i}=0$ and subtracting one from the largest remaining $d_{i}$ degrees in $d$.  This is equivalent to connecting vertex $u=\IV{i}$ to the vertices of largest degree, not choosing $u$ to avoid self-loops.  Note that there may be more than $d_{i}$ degrees to choose from the largest degrees.  To make this more precise, let $K$ contain $d_{i}$ indices, not including $i$, of maximal degree in $d$, calling $[i,K]$ a {\it maximal index pair}.  
We can represent the step above by
\[\hat{d} = \HH(d,i,K) \equiv d - \sum_{j\in K} e^{ij},\]
with $e^{ij} \equiv \{\delta_{ik}+\delta_{jk}\}_{k=1}^{N}$ and $\delta_{ij}$ the Kronecker delta operator.  We will call the operator $\Delta$ the Havel-Hakimi operator, in reference to the following well-known theorem:
\begin{theorem}[Havel-Hakimi \cite{Havel:1955p7653, hakimi:496, hakimi:135}]
\label{thm:HH}
If $d$ is an integer sequence, then for any maximal index pair $[i,K]$, $d$ is graphic $\Longleftrightarrow$ $\HH(d,i,K)$ is graphic.
\end{theorem}

We now address integer-pair sequences.  In Theorem~\ref{thm:erdos}, we needed the integer sequence to be non-increasing.  To determine when an integer-pair sequence is digraphic, we will need the integer-pair sequence to be non-increasing relative to the lexicographical ordering.
\begin{definition}
An integer-pair sequence $d = \{(d^{+}_{i}, d^{-}_{i})\}_{i=1}^{N}$ is non-increasing relative to the {\bf positive lexicographical ordering} if and only if $d^{+}_{i} \geq d^{+}_{i+1}$, with $d^{-}_{i} \geq d^{-}_{i+1}$ when $d^{+}_{i} = d^{+}_{i+1}$.  In this case, we will call $d$ {\bf positively ordered} and denote the ordering by $d_{i} \pgeq d_{i+1}$.  We say $d$ is non-increasing relative to the {\bf negative lexicographical ordering} by giving preference to the second index, calling $d$ in this case {\bf negatively ordered} and denoting the ordering by $d_{i} \ngeq d_{i+1}$.
\end{definition}

For a given integer-pair sequence $d = \{(d^{+}_{i},d^{-}_{i})\}_{i=1}^{N}$, define the sequences $\bar{d} = \{(\bar{d}^{+}_{i}, \bar{d}^{-}_{i})\}_{i=1}^{N}$ and $\ubar{d} = \{(\ubar{d}^{+}_{i}, \ubar{d}^{-}_{i})\}_{i=1}^{N}$ to be the positive and negative orderings of $d$, respectively.

We have the following theorem by Fulkerson:
\begin{theorem}[Fulkerson \cite{Fulkerson:1960p3387}]
\label{thm:fulkerson}
An integer-pair sequence $d$ is digraphic if and only if $\sum_{i=1}^{N} d^{+}_{i} = \sum_{i=1}^{N} d^{-}_{i}$ and for $k=1,\ldots,N$,
\begin{gather*}
\sum_{i=1}^{k} \min[\bar{d}^{-}_{i},k-1] + 
\sum_{i=k+1}^{N}\min[\bar{d}^{-}_{i},k] \geq \sum_{i=1}^{k}\bar{d}^{+}_{i}\\
(or)\\
\sum_{i=1}^{k} \min[\ubar{d}^{+}_{i},k-1] + \sum_{i=k+1}^{N}\min[\ubar{d}^{+}_{i},k] \geq \sum_{i=1}^{k}\ubar{d}^{-}_{i}.
\end{gather*}
\end{theorem}
Note that the last inequality for $k=N$ is actually an equality by the condition $\sum_{i=1}^{N} d^{+}_{i} = \sum_{i=1}^{N} d^{-}_{i}$.

We will actually use a restatement of Theorem~\ref{thm:fulkerson}.  For that restatement, we need some definitions.  Given an integer sequence $a$, define the {\it corrected conjugate sequence} $a^{\prime\prime}$ \cite{Berge:1973yq} by
\[
a_{k}^{\prime\prime} = |I_{k}| + |J_{k}|,
\]
where
\begin{eqnarray*}
I_{k} & = & \{i \ | \ i < k \ \ \mathrm{and} \ \ a_{i} \geq k-1\}, \\
J_{k} & = & \{i \ | \ i > k \ \ \mathrm{and} \ \ a_{i} \geq k\}.
\end{eqnarray*}
The numbers $a^{\prime\prime}$ can be represented by what is known as the {\it corrected Ferrers diagram} \cite{Berge:1973yq}, shown with the example sequence $a = (4 \ 3 \ 4 \ 2 \ 1)$ below.
\[
\begin{array}{c|ccccc|c}
      & a_{1}^{\prime\prime} & a_{2}^{\prime\prime} & a_{3}^{\prime\prime} & a_{4}^{\prime\prime} & a_{5}^{\prime\prime} \\
      \hline
a_{1} & \circ   & \bullet & \bullet & \bullet & \bullet & 4 \\
a_{2} & \bullet & \circ   & \bullet & \bullet & \circ   & 3 \\
a_{3} & \bullet & \bullet & \circ   & \bullet & \bullet & 4 \\
a_{4} & \bullet & \bullet & \circ   & \circ   & \circ   & 2 \\
a_{5} & \bullet & \circ   & \circ   & \circ   & \circ   & 1 \\ \hline
      & 4 & 3 & 2 & 3 & 2
\end{array}
\]
If in the $i$-th row $a_{i}$ solid dots are filled in from left to right, making sure we skip the $i$-th column, then the value $a_{k}^{\prime\prime}$ is found by simply counting the number of solid dots in the $k$-th column, giving $a^{\prime\prime} = (4 \ 3 \ 2 \ 3 \ 2)$.
The partial sums of the corrected conjugate sequence are given by
\begin{equation}
\label{equ:corrconj}
\sum_{i=1}^{k}a^{\prime\prime}_{i} = \sum_{i=1}^{k} \min[a_{i},k-1] + \sum_{i=k+1}^{N}\min[a_{i},k].
\end{equation}
For an integer-pair sequence $d$, define the {\it slack} sequences $\bar{s}$ and $\ubar{s}$ by
\begin{eqnarray*}
\bar{s}_{k} & = & \sum_{i=1}^{k}[\bar{d}^{-}]_{i}^{\prime\prime}-\sum_{i=1}^{k}\bar{d}^{+}_{i} \ \ \mbox{with} \ \ \bar{s}_{0} \equiv 0, \\
\ubar{s}_{k} & = & \sum_{i=1}^{k}[\ubar{d}^{+}]_{i}^{\prime\prime}-\sum_{i=1}^{k}\ubar{d}^{-}_{i} \ \ \mbox{with} \ \ \ubar{s}_{0} \equiv 0.
\end{eqnarray*}
Note that it is possible for the slack sequences to be negative.  However, if the integer-pair sequence is digraphic, then they will both be non-negative, as can be seen by the following simple restatement of Theorem~\ref{thm:fulkerson} using the slack sequence and Eq.~(\ref{equ:corrconj}):
\begin{cor}
An integer-pair sequence $d$ is digraphic if and only if $\sum_{i=1}^{N} d^{+}_{i} = \sum_{i=1}^{N} d^{-}_{i}$ and for $k=1,\ldots,N$,
\begin{gather*}
\sum_{i=1}^{k}[\ubar{d}^{+}]_{i}^{\prime\prime} \geq \sum_{i=1}^{k}\ubar{d}_{i}^{-} \ \ \Leftrightarrow \ \ \ubar{s}_{k} \geq 0, \\
(or)\\
\sum_{i=1}^{k}[\bar{d}^{-}]_{i}^{\prime\prime} \geq \sum_{i=1}^{k}\bar{d}_{i}^{+} \ \ \Leftrightarrow \ \ \bar{s}_{k} \geq 0.
\end{gather*}
\end{cor}
As in Theorem~\ref{thm:fulkerson}, we have equality for $k=N$, showing $\ubar{s}_{N}=\bar{s}_{N} = 0$.

We also have a theorem by Kleitman and Wang \cite{Kleitman:1973p4546} for constructing simple digraphs from an integer-pair sequence.  In this proof, and in many proofs in this paper, there are some common arc switches that are used throughout, notably 2- and 3-switches.  These are defined in full generality in the following definition \cite{McDonald:2007p6830}.
\begin{definition}
Let $\DG=(V,\DE)$ be a digraph with $\{x_{1},\ldots,x_{n}\}\subseteq V$ and $\{y_{1},\ldots,y_{n}\}\subseteq V$ two sets of $n$ distinct vertices such that  $\DS=\{(x_{i},y_{i})\}_{i=1}^{n}\subseteq \DE$.  Given the reordering operator $\pi(\DS) \equiv \{(x_{i},y_{\mod(i,n)+1})\}_{i=1}^{n}$, if we include the constraints

(i) $x_{i} \neq y_{i} \neq y_{\mod(i,n)+1}$

(ii) $\pi(\DS) \subseteq \DE^{C}$

\noindent then we can define the {\bf $\bm n$-switch} operator $\sigma_{n}$ by
\[\sigma_{n}(\DG,\DS) = \DG^{\prime} = (V,\DE^{\prime})\]
with
\[\DE^{\prime} = \left(\DE\setminus \DS\right)\cup \pi(\DS).\]
\end{definition}


We will sometimes denote it by $\sigma_{n}(S)$, with the graph $\DG$ understood by context.  Note that an $n$-switch preserves degrees, with the constraints $(i)$ and $(ii)$ imposed to avoid self-loops and multi-arcs, respectively.  As examples, a 2-switch $\sigma_{2}\bigl((x_{1},y_{1}),(x_{2},y_{2})\bigr)$ is shown by the following diagram
\begin{center}
\medskip
\includegraphics{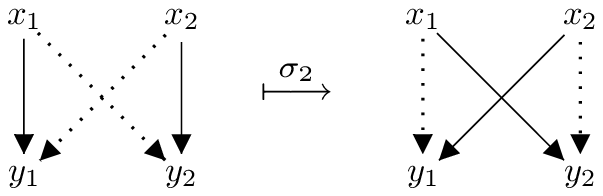}
\medskip
\end{center}
while a 3-switch $\sigma_{3}\bigl((x_{1},y_{1}),(x_{2},y_{2}),(x_{3},y_{3})\bigr)$ is given by
\begin{center}
\medskip
\includegraphics{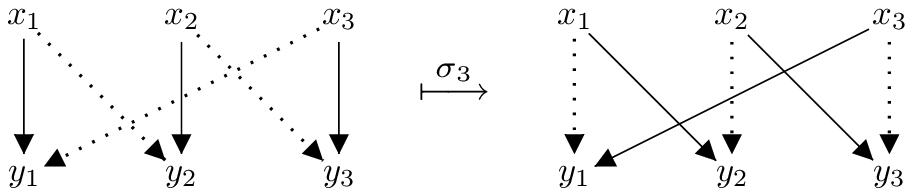}
\medskip
\end{center}

%
%

\begin{definition}

(i) Given an index $i$ of $d$, $\sK = \{\Kp, \Km\}$ are {\bf maximal index sets} if and only if $\Kp$ ($\Km$) contains $d_{i}^{-}$ ($d_{i}^{+}$) indices, not including $i$, of maximal degree relative to the positive (negative) ordering.  We will call $[i,\sK]$ a {\bf maximal index pair}.

(ii) Given a vertex $u$ of $G_{d}$, $\sM = \{\Mp, \Mm\}$ are {\bf maximal vertex sets} if and only if $\L(\sM) \equiv \{\L(\Mp),\L(\Mm)\}$ are maximal index sets for $\L(u)$.  We will call $[u,\sM]$ a {\bf maximal vertex pair}.
\end{definition}

Let $d = \{(d^{+}_{i}, d^{-}_{i})\}_{i=1}^{N}$ be a non-trivial integer-pair sequence of length $N \geq 2$.  Choose a maximal index pair $[i,\sK]$ of $d$ with $d^{+}_{i} > 0$.  Construct $\hat{d}$ by setting $d^{+}_{i} = 0$ and subtracting 1 from all degrees with indices in the set $\Km$.  
We can represent the step above by
\[\hat{d} = \HH^{+}(d,i,\sK) \equiv d - \sum_{j\in\Km}\xi^{ij},\]
with $\xi^{ij} \equiv \{(\delta_{ik},\delta_{jk})\}_{k=1}^{N}$.  We include a modified proof of Kleitman and Wang's theorem.
\begin{theorem}[Kleitman and Wang \cite{Kleitman:1973p4546}]
If $d$ is an integer-pair sequence, then for any maximal index pair $[i,\sK]$, $d$ is digraphic $\Longleftrightarrow$ $\HH^{+}(d,i,\sK)$ is digraphic.  
\label{thm:KW}
\end{theorem}

\bnprf
If $\hat{d}$ is digraphic, then given a digraph $\DG_{\hat{d}}$, construct a digraph $\DG_{d}$ by adding a vertex $v$ to $\DG_{\hat{d}}$ and connecting $v$ to the vertices of highest degree in $\DG_{\hat{d}}$.  This gives the degree sequence $d$ with realization $\DG_{d}$, showing $d$ is digraphic.

Now suppose $d$ is digraphic and denote a realization of $d$ by $\DG_{d} = (V,\DE)$.  Let $[i,\sK]$ be a maximal index pair of $d$ such that $n = d_{i}^{+} > 0$,  with $[u,\sM]=[\IV{i},\IV{\sK}]$ the corresponding maximal vertex pair.  If we can switch {\de}s of $\DG_{d}$ such that $(u,\Mm) \subseteq \DE$, then we can remove those {\de}s, thereby realizing a digraph $\DG_{\hat{d}}$.

If $(u,\Mm) \subseteq \DE$, then we're done, so suppose there is a $w \in \Mm$ such that $(u,w)\notin \DE$.  Let $v \notin \Mm$ such that $(u,v)\in \DE$.  Since $v \notin \Mm$, $d_{v}^{-} \leq d_{w}^{-}$.  If $d_{v}^{-} < d_{w}^{-}$, then there is a vertex $x \neq v$ such that $(x,w)\in \DE$ and $(x,v)\notin \DE$.  Without affecting degrees, we can perform the 2-switch $\sigma_{2}\bigl((u,v),(x,w)\bigr)$.  We are left with the case $d_{v}^{-} = d_{w}^{-}$.

Since $(u,v)\in \DE$, $d_{v}^{-} = d_{w}^{-} \geq 1$.  So if $(v,w)\notin \DE$, there is an $x$ such that $(x,w)\in \DE$ and $(x,v)\notin \DE$.  As before, we can perform the 2-switch $\sigma_{2}\bigl((u,v),(x,w)\bigr)$.

Now consider the case $(v,w)\in \DE$.  We must move to $d^{+}$ to give us information, and so we use the fact that $d_{v} \nleq d_{w}$.  From this and $d_{v}^{-} = d_{w}^{-}$, $d_{v}^{+} \leq d_{w}^{+}$.  Since $(v,w)\in \DE$, $d_{v}^{+} \geq 1$, and so $d_{w}^{+} \geq 1$.  If $(w,v)$ is an {\de }, then $d_{v}^{-} \geq 2$ and thus $d_{w}^{-} \geq 2$.  Thus, there exists an $x \neq u, v$ such that $(x,v)\notin \DE$ and $(x,w)\in \DE$.  We thus perform the 2-switch $\sigma_{2}\bigl((u,v),(x,w)\bigr)$.

Assuming there is not an {\de } $(w,v)$, we consider the case where $(v,u)\in \DE$.  Since $d_{v}^{+} \leq d_{w}^{+}$, we have $d_{v}^{+} \geq 2$ and thus $d_{w}^{+} \geq 2$.  With $(w,v)\notin \DE$, there is an $x \neq u, v$ such that $(v,x)\notin \DE$ and $(w,x)\in \DE$.  We thus perform the 3-switch $\sigma_{3}\bigl((u,v),(v,w),(w,x)\bigr)$.

We are left with existence of the directed path $(u,v,w)$ with $\{(w,v),(v,u)\} \subseteq \DE^{C}$.   Again, since $d_{v}^{+} \leq d_{w}^{+}$ and $(v,w)\in \DE$, there is a vertex $x \neq v$ such that $(w,x)\in \DE$ and $(v,x)\notin \DE$, giving the path $(u,v,w,x)$.  We then perform the 3-switch $\sigma_{3}\bigl((u,v),(v,w),(w,x)\bigr)$.

This covers all the cases, and we see that we finally have an {\de } $(u,w)\in\DE$.  Repeating for the remaining vertices in $\Mm$ and removing all {\de}s $(u,\Mm) \subseteq \DE$, we have a realization of $\DG_{\hat{d}}$.
\edprf

We can also construct the residual degree sequence $\hat{d}$ by setting $d^{-}_{i} = 0$ and subtracting 1 from the degrees corresponding to the index set $\Kp$.
We then represent the Havel-Hakimi step by
\[\hat{d} = \HH^{-}(d,i,\sK) \equiv d - \sum_{j\in\Kp}\xi^{ji}.\]
A similar proof above shows that $d$ is digraphic if and only if $\hat{d}$ is digraphic.  For a given index $i$, we can do either $\HH^{-}$ or $\HH^{+}$, both of which may give different realizations.  Also, for integer sequences, we actually removed the $i$-th index after one Havel-Hakimi step, which may not happen for integer-pair sequences in Kleitman and Wang's algorithm.  There are two natural ideas to try so that after one ``step'', we remove the $i$-th index.  The first is to count one ``step'' as two Havel-Hakimi steps in serial, i.e. $\HH^{+} \circ \HH^{-}$ or $\HH^{-} \circ \HH^{+}$.  As long as you use the same index $i$ at each step, the new {\it serial} Havel-Hakimi operator will remove the $i$-th index after one step.  The second idea is somewhat more natural (why it is more natural will be explained below) and consists of a {\it parallel} Havel-Hakimi step $\HH^{\pm}$.  Given an index pair $[i,\sK]$, 
define the parallel Havel-Hakimi step by
\[\hat{d} = \HH^{\pm}(d,i,\sK) \equiv d - \sum_{j\in\Kp}\xi^{ji} - \sum_{j\in\Km}\xi^{ij}.\]
The reason this is considered more natural is shown by the example of an integer sequence $d = (1\,1\,1\,1)$ and its extension to an integer-pair sequence
\[
\tilde{d} = 
\left( \begin{array}{cccc}
           1 & 1 & 1 & 1 \\
           1 & 1 & 1 & 1
       \end{array}
\right).
\]

The parallel Havel-Hakimi operator acting on $\tilde{d}$ will give the same realization as Havel-Hakimi acting on $d$ when bidirectional {\de}s are identified with undirected edges and $\Kp=\Km$.  This is not the case for the serial Havel-Hakimi operator, whose realization is shown below as the digraph on the right, with the realization for the parallel operator on the left.
\medskip
\begin{center}
\includegraphics{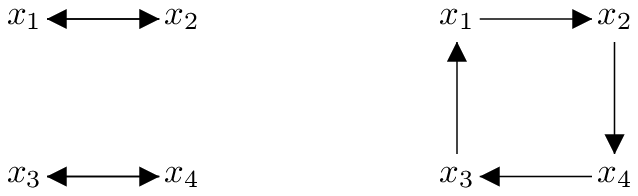}
\end{center}
\medskip
It is clear that the serial Havel-Hakimi operator has a theorem analogous to Kleitman and Wang's, since it is just a composition of two Havel-Hakimi operators, both of which satisfy Theorem~\ref{thm:KW}.  Is there such a theorem for the parallel operator?  The answer is no as it stands.  When and how does it break down?  It is clear that if $\hat{d}$ is digraphic, then $d$ is digraphic, so the problem is when $d$ is digraphic but $\hat{d}$ is not digraphic.  

The simplest example of such an ill-defined degree sequence is
\[ \left( \begin{array}{ccc}
               1 & 1 & 1 \\
               1 & 1 & 1
            \end{array}
   \right),
\]
which is the degree sequence for a directed 3-cycle.  If we consider the first index, then $\hat{d} = \HH^{\pm}(d,1,\sK)$ will be digraphic only if $\Kp \cap \Km = \emptyset$.  Otherwise,
$\hat{d}$ will be, for example, 
\[ \left( \begin{array}{ccc}
               0 & 0 & 1 \\
               0 & 0 & 1
            \end{array}
   \right)
\]
which is not digraphic.  Notice that all of the indices are problematic.  Another example, shown to illustrate that these ill-defined degree sequences are non-trivial, is the degree sequence
\[ \left( \begin{array}{cccccc}
              1 & 2 & 0 & 2 & 5 & 2 \\
              2 & 3 & 1 & 3 & 0 & 3
          \end{array}
   \right).
\]
The corresponding labeled vertex set is given by $\{x_{1},x_{2},x_{3},x_{4},x_{5},x_{6}\}$ with a digraph realization given by
\begin{center}
\includegraphics{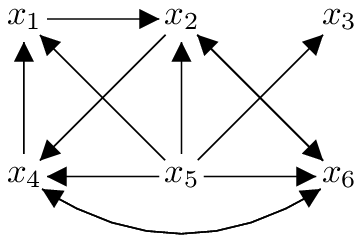}
\end{center}
In this example, the first index is the only ill-defined index for $\HH^{\pm}$.  In other words, for indices $i\in\{2,\ldots,6\}$, $\HH^{\pm}(d,i,\sK)$ is digraphic for {\it any} choice of maximal index sets $\sK$.  This leads to a definition:

\begin{definition}

(i) A digraphic degree sequence $d \in \gD$ if and only if there exists an {\bf ill-defined} maximal index pair $[i,\sK]$ such that $\HH^{\pm}(d,i,\sK)$ is not digraphic.

(ii) A digraph $\DG \in \gR$ if and only if there exists an {\bf ill-defined} maximal vertex pair $[u,\sM]$ such that for all graphs $\DG^{\prime} \in R(d_{\DG})$, $(\Mp,u,\Mm) \not\subseteq \DP(\DG^{\prime})$.
\label{def:gDgR}
\end{definition}

\begin{lem}
Definitions~\ref{def:gDgR}(i) and (ii) are equivalent, i.e.
\[\gR = R(\gD) \ \ \mbox{and} \ \ \gD = d_{\gR}\]
with an index pair ill-defined if and only if the corresponding vertex pair is ill-defined.
\label{lem:gDgR}
\end{lem}
\bnprf
We will prove the first equality since the second will follow directly.  Let $d \in \gD$.  We have an ill-defined index pair $[i,\sK]$ such that $\HH^{\pm}(d,i,\sK)$ is not digraphic.  However, if $[u,\sM] = [\IV{i},\IV{\sK}]$, then for every digraph $\DG \in R(d)$, $(\Mp,u,\Mm) \not\subseteq \DP(\DG)$, for otherwise we would have a realization of $\HH^{\pm}(d,i,\sK)$ by removing those {\de}s from $\DG$.  Thus, $[u,\sM]$ is an ill-defined vertex pair and hence $R(\gD) \subseteq \gR$.

Suppose now $\DG \in \gR$ and $d_{\DG} \notin \gD$.  Let $[u,\sM]$ be an ill-defined vertex pair and $[i,\sK]$ be a corresponding maximal index pair.  Since $d_{\DG} \notin \gD$, $\hat{d} = \HH^{\pm}(d_{\DG},i,\sK)$ is digraphic, and so let $\DG^{\prime} \in R(\hat{d})$.  If we let $\DG^{\prime\prime} = (V^{\prime\prime},\DE^{\prime\prime})$ be such that $V^{\prime\prime} = V(\DG^{\prime})\cup \{u\}$ and $\DE^{\prime\prime} = \DE(\DG^{\prime})\cup (\Mp,u)\cup (u,\Mm)$, then $\DG^{\prime\prime} \in R(d_{G})$ and $(\Mp,u,\Mm) \subseteq \DP(\DG^{\prime\prime})$, which is a contradiction.  Thus, $[i,\sK]$ is an ill-defined index pair and $\gR \subseteq R(\gD)$.

\edprf

Before moving on to the main theorems of the paper which characterize the ill-defined degree sequences and their digraph realizations, we prove in Theorem~\ref{thm:digraphic} that given any index $i$ of $d$, we can choose a particular maximal index set $\sK$ so that $[i,\sK]$ is well-defined.  For this we need two results, the first of which shows that the difficulty in the parallel Havel-Hakimi algorithm lies in the intersection $\Mp\cap\Mm$ of the maximal sets.
\begin{lem}
If $d \in \gD$, then for every ill-defined maximal index pair $[i,\sK]$ with corresponding maximal vertex pair $[u,\sM]$, there exists a digraph $\DG\in R(d)$ with vertex set $C=\{u,v,w\}$ such that

(i) $(\Mp,u,\Mm\setminus \Mp) \subseteq \DP(\DG)$ and

(ii) $\DG[C] \cong \CC$ with $w \in \Mp\cap\Mm$, $v\notin\Mp\cup\Mm$, and $d^{-}_{v}=d^{-}_{w}$.

\label{lem:arbcon}
\end{lem}

\bnprf
Let $\DG = (V,\DE)\in \gR$.  We can apply the techniques from Theorem \ref{thm:KW} to the set $\Mp$ to arrive at a digraph $\DG^{\prime}(V,\DE^{\prime})\in R(d_{\DG})$ such that $(\Mp,u) \subseteq \DE^{\prime}$.  We apply the techniques from Theorem \ref{thm:KW} again to the set $\Mm\setminus \Mp$, except we have to be careful along the way in the proof to make sure we don't disconnect any {\de}s from $\Mp$ to $u$.  Let $w^{\prime} \in \Mm\setminus \Mp$.  Upon reviewing the proof, we see that one of the situations where we need to be careful is when $(v^{\prime},u)$ is an {\de }, which implies $v^{\prime} \in \Mp$ since we've already made that connection.  However, it is immediately clear that when we perform the 3-switch in the proof, the {\de } $(v^{\prime},u)$ is preserved.  The second case that could pose trouble is when $(w^{\prime},u) \in \DE^{\prime}$, i.e. when we have a directed 3-cycle.  But again, this means $w^{\prime} \in \Mp$, which can't happen by assumption.  This covers all cases, showing we can successfully find a digraph $\DG^{\prime\prime}(V,\DE^{\prime\prime})\in R(d_{\DG})$ through a series of arc switches so that $(i)$ is satisfied, i.e.
\[(\Mp,u,\Mm \setminus \Mp) \subseteq \DP(\DG^{\prime\prime}).\]
For $(ii)$, since $d \in \gD$, we must have $w \in \Mp\cap\Mm$ with $(w,u)\in\DE$, which is again the case where there is a $v\notin\Mp\cap\Mm$ such that $C=\{u,v,w\}$ induces a directed 3-cycle with $d_{v}^{-} = d_{w}^{-}$.
\edprf





\begin{lem}
If $d\in \gD$ with an ill-defined maximal index pair $[i,\sK]$ and corresponding maximal vertex pair $[u,\sM]$, then there exists a $\DG\in R(d)$ that satisfies Lemma~\ref{lem:arbcon} with the added property that for $x \in V \setminus C$,

(i) \ $(x,w)\in \DE \ \Longleftrightarrow \ (x,v)\in \DE$,

(ii) $(w,x)\in \DE \ \Longleftrightarrow \ (v,x)\in \DE$.

\label{lem:vweqz}
\end{lem}

\bnprf
$(i)$ If we have $(x,w)\in \DE$ but $(x,v)\notin \DE$, then we can perform the 2-switch $\sigma_{2}\bigl((u,v),(x,w)\bigr)$.  Now suppose $(x,v)\in \DE$ but $(x,w)\notin \DE$.  Since $d_{v}^{-} = d_{w}^{-}$, there is a $z$ such that $(z,w)\in \DE$ and $(z,v)\notin \DE$.  We can thus perform the 2-switch $\sigma_{2}\bigl((u,v),(z,w)\bigr)$.

$(ii)$ If $(w,x)\in \DE$ but $(v,x)\notin \DE$, then we perform the 3-switch $\sigma_{3}\bigl((u,v),(v,w),(w,x)\bigr)$.  Conversely, supposing $(v,x)\in \DE$ and $(w,x)\notin \DE$, since $d_{v}^{+} \leq d_{w}^{+}$, there is a $z$ such that $(w,z)\in \DE$ and $(v,z)\notin \DE$.  Similar to above, we perform the 3-switch $\sigma_{3}\bigl((u,v),(v,w),(w,z)\bigr)$.
\edprf


\begin{theorem} 
\label{thm:digraphic}
For every index $i$ of a digraphic degree sequence $d$, there exists maximal index sets $\sK$ for $i$ such that $\HH^{\pm}(d,i,\sK)$ is digraphic.
\end{theorem}

\bnprf
If $d \notin \gD$, then the statement is trivial by definition, so suppose $d \in \gD$.  Let $[i,\sK]$ be an ill-defined maximal index pair such that $\HH^{\pm}(d,i,\sK)$ is not digraphic.  Let $[u,\sM] = [\IV{i},\IV{\sK}]$.  By Lemma~\ref{lem:arbcon}, there is a $\DG\in R(d)$ such that $(\Mp,u,\Mm\setminus \Mp)\subseteq \DP(\DG)$ with vertices $\{v,w\}\subseteq V(\DG)$ such that $C=\{u,v,w\}$ is $\CC$, $w \in \Mp\cap \Mm$ and $v \notin \Mp\cup \Mm$.  At this point, we have $d_{v}^{-} = d_{w}^{-}$ and $d_{v}^{+} \leq d_{w}^{+}$.  But by Lemma~\ref{lem:vweqz}, we have the stronger statement
\[d_{v} = d_{w}.\]
Let $\hat{\sM} \equiv \{\hMp,\hMm\}$ be such that
\begin{eqnarray*}
\hMp & = & \Mp \\
\hMm & = & (\Mm \setminus \{w\}) \cup \{v\}.
\end{eqnarray*}
Repeating for the remaining vertices in $\hMp\cap\hMm$ and letting $\tilde{\sM}$ be the resulting maximal vertex sets, with $\tilde{\sK} = \VI{\tilde{\sM}}$, we have $\tilde{K}^{+}\cap\tilde{K}^{-} = \emptyset$ and $\HH^{\pm}(d,i,\tilde{\sK})$ digraphic.
\edprf



\section{Characterizations of ill-defined degree sequences and their realizations} \label{sec:characterizations}

Ultimately, to make the parallel Havel-Hakimi algorithm well-defined we need to be able to detect the ill-defined degree sequences $\gD$, and in particular the ill-defined indices.  It is not enough, however, to just know the degree sequences that cause problems for the algorithm, since once they are detected we need to know how to make connections in the resulting digraph realizations $\gR$ in order to proceed.  This section accomplishes both of these goals by proving two main theorems (Theorems~\ref{thm:charequ} and \ref{thm:equiv}) which provide a degree sequence characterization for $\gD$, as well as a structural characterization of all digraph realizations $\gR$.

\begin{figure}
\scalebox{0.65}{
\begin{picture}(80,240)(-130,0)
\includegraphics{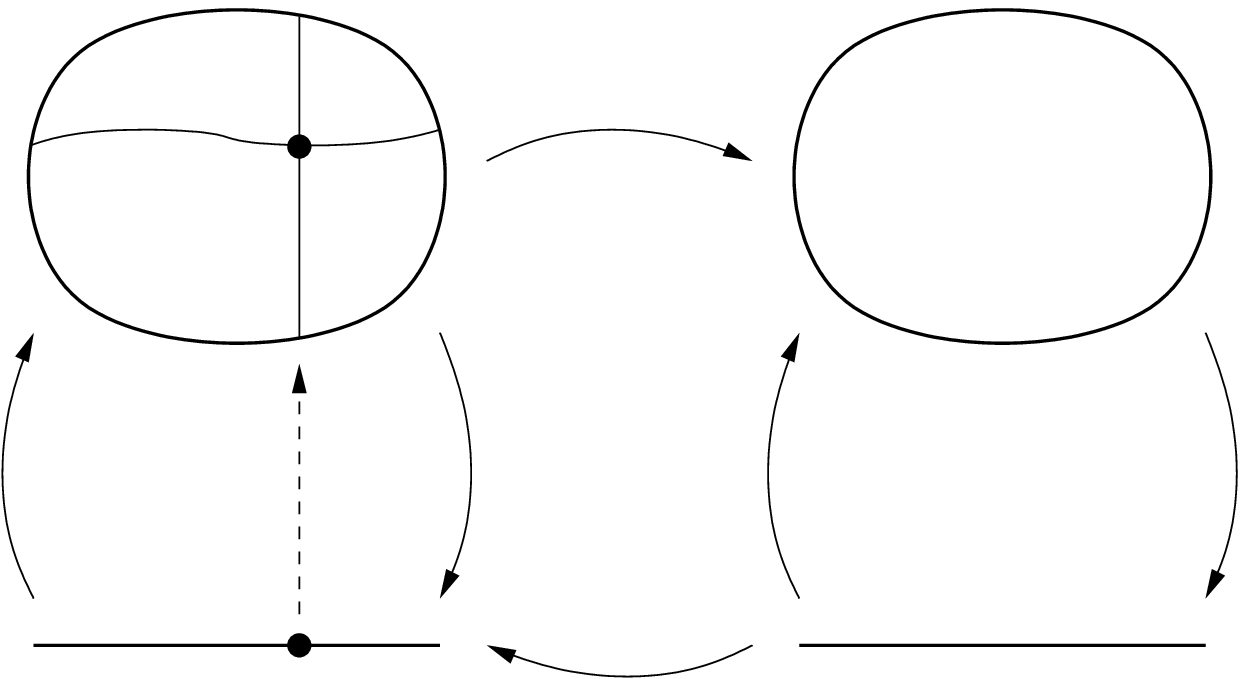}

\put(-340,225){\scalebox{1.5}{Ill-defined sets}}
\put(-127,225){\scalebox{1.5}{Characterizations}}
\put(-470,145){\scalebox{1.5}{Digraphs}}
\put(-470,5){\scalebox{1.5}{$\stackrel{\textstyle \mathrm{Degree}}{\textstyle \mathrm{sequences}}$}}

\put(-384,177){\scalebox{2}{$\stackrel{\textstyle \gR}{\scriptscriptstyle [\ref{def:gDgR}]}$}}
\put(-5,177){\scalebox{2}{$\stackrel{\textstyle \gT}{\scriptscriptstyle [\ref{def:gT}]}$}}
\put(-5,-15){\scalebox{2}{$\stackrel{\textstyle \gS}{\scriptscriptstyle [\ref{def:gS}]}$}}
\put(-384,-15){\scalebox{2}{$\stackrel{\textstyle \gD}{\scriptscriptstyle [\ref{def:gDgR}]}$}}
\put(-325,137){\scalebox{1.5}{$\stackrel{\textstyle \gQ}{\scriptscriptstyle [\ref{def:gQ}]}$}}
\put(-270,125){\scalebox{1.5}{$R(d)$}}
\put(-278,-10){\scalebox{1.5}{$d$}}

\put(-203,165){\scalebox{1.5}{$\stackrel{\textstyle \gQ\subseteq\gT}{\scriptscriptstyle [\ref{lem:gQgT}]}$}}
\put(-200,-23){\scalebox{1.5}{$\stackrel{\textstyle \gS\subseteq\gD}{\scriptscriptstyle [\ref{lem:gSgR}]}$}}

\put(-394,50){\scalebox{1.5}{\colorbox{white}{$\stackrel{\textstyle 
R(\gD)=\gR}{\scriptscriptstyle [\ref{lem:gDgR}]}$}}}
\put(-251,50){\scalebox{1.5}{\colorbox{white}{$\stackrel{\textstyle d_{\gR}=\gD}{\scriptscriptstyle [\ref{lem:gDgR}]}$}}}

\put(-173,50){\scalebox{1.5}{\colorbox{white}{$\stackrel{\textstyle 
R(\gS)\subseteq\gT}{\scriptscriptstyle [\ref{lem:gRgT}]}$}}}
\put(-30,50){\scalebox{1.5}{\colorbox{white}{$\stackrel{\textstyle d_{\gT}\subseteq\gS}{\scriptscriptstyle [\ref{lem:gTgS}]}$}}}
\end{picture}
}
\vspace{0.15in}
\caption{Diagram showing relationships between the 5 sets $\gD$, $\gR$, $\gQ$, $\gT$, and $\gS$, with references to the corresponding definitions and theorems in brackets.}
\label{fig:diagram}
\end{figure}

The proof technique involves defining three sets: the set of degree sequences $\gS$ defined in Definition~\ref{def:gS} which we show in Theorem~\ref{thm:equiv} is equivalent to the ill-defined sequences $\gD$; the set of digraphs $\gT$ defined in Definition~\ref{def:gT} which we show in Theorem~\ref{thm:equiv} is equivalent to the ill-defined digraphs $\gR$; and an intermediate set of digraphs $\gQ \subseteq \gR$.  The relationships between $\gS$, $\gQ$, $\gT$, $\gD$ and $\gR$ are shown in Figure~\ref{fig:diagram}, giving a graphical outline with references to the definitions and lemmas necessary to prove in Theorem~\ref{thm:equiv} that $\gS=\gD$ and $\gT=\gR$.  The subsections are as follows: Section \ref{sec:gQgT} defines the set of digraphs $\gQ$ and $\gT$, showing that by definition $\gQ \subseteq \gR$ and proving $\gQ \subseteq \gT$; Section \ref{sec:gSgR} defines the set of degree sequences $\gS$ and shows $\gS \subseteq \gD$; and Section \ref{sec:gTgS} proves that $d_{\gT} \subseteq \gS$ and $R(\gS) \subseteq \gT$.  The results of these three subsections allow us to prove the following main theorems of the paper.
\begin{theorem}
\[R(\gS) = \gT \ \ \mbox{and} \ \ d_{\gT} = \gS\]
\label{thm:charequ}
\end{theorem}

\bnprf
By Lemma~\ref{lem:gTgS} and \ref{lem:gRgT},
\[
R(\gS) \subseteq \gT \Longrightarrow \gS \subseteq d_{\gT} \subseteq \gS
\]
and thus $d_{\gT} = \gS$.  We also have
\[
\gT \subseteq R(d_{\gT}) = R(\gS) \subseteq \gT
\]
and thus $R(\gS) = \gT$.
\edprf

Finally, we can easily prove equivalence of the ill-defined sets and the characterizations as follows:

\begin{theorem}
\[\gD = \gS \ \ \mbox{and} \ \ \gR = \gT\]
\label{thm:equiv}
\end{theorem}

\bnprf
By the definition of $\gQ$, Lemmas~\ref{lem:gQgT} and \ref{lem:gSgR}, and Theorem~\ref{thm:charequ}, we have
\[
\gD = d_{\gQ} \subseteq d_{\gT} = \gS \subseteq \gD
\]
and thus $\gS = \gD$.  We also have by Lemma~\ref{lem:gDgR} and Theorem~\ref{thm:charequ}
\[
\gR = R(\gD) = R(\gS) = \gT.
\]
\edprf


%


\subsection{Structural characterization}\label{sec:gQgT}

In this section we define the set $\gT$ which we show in Theorem~\ref{thm:equiv} is the set of digraphs that constitute the structural charactarization of the ill-defined digraphs $\gR$.  We will first show there is a set $\gQ \subseteq \gR$ whose digraphs have a structural characterization by a digraph decomposition using $M$-partitions \cite{Feder:2003p7761, Cameron:2007p7831}, and then show in Lemma~\ref{lem:gQgT} that $\gQ\subseteq\gT$.  An $M$-partition of a digraph $\DG$ is a partition of the vertex-set $V(\DG)$ into $k$ disjoint classes $\{X_{1},\ldots,X_{k}\}$, where the {\de } constraints within and between classes are given by a $k\times k$ matrix $M$ with elements in $\{0,1,*\}$.  $M_{ii}$ equal to $0$, $1$, or $*$ corresponds to $\DG[X_{i}]$ being an independent set, clique, or arbitrary subgraph, respectively.  Similarly, for $i\neq j$, $M_{ij}$ equal to $0$, $1$, or $*$ corresponds to $\DG[X_{i},X_{j}]$ having no {\de}s from $X_{i}$ to $X_{j}$, all {\de}s from $X_{i}$ to $X_{j}$, and no constraint on {\de}s from $X_{i}$ to $X_{j}$, respectively.

\begin{definition}
Let $\DG=(V,\DE)$ be a digraph with vertex set $C=\{u,v,w\}\subseteq V$ such that $\DG[C]\cong \CC$ and an $M$-partition of $\DG[V\setminus C]$ with vertex classes given by $\{\I,\U,\Z,\Cp,\Cm,\C\}$.  With $Z\equiv\{v,w\}$, each class defines how its elements relate to $C$ as follows:
\begin{eqnarray*}
	\I \ & \equiv & \ \{ x \in V\setminus C \ : \ (x,C)\cup(C,x) \subseteq \DE^{C} \} \\
	\U \ & \equiv & \ \{ x \in V\setminus C \ : \ \{(x,u),(u,x)\} \subseteq \DE \ \mbox{and} \ (x,Z)\cup(Z,x) \subseteq \DE^{C} \} \\
	\Z \ & \equiv & \ \{ x \in V\setminus C \ : \ (x,Z)\cup(Z,x) \subseteq \DE \ \mbox{and} \ \{(x,u),(u,x)\} \subseteq \DE^{C} \} \\
	\Cp \ & \equiv & \ \{ x \in V\setminus C \ : \ (x,C) \subseteq \DE \ \mbox{and} \ (C,x) \subseteq \DE^{C} \} \\
	\Cm \ & \equiv & \ \{ x \in V\setminus C \ : \ (C,x) \subseteq \DE \ \mbox{and} \ (x,C) \subseteq \DE^{C} \} \\
	\C \ & \equiv & \ \{ x \in V\setminus C \ : \ (x,C)\cup(C,x) \subseteq \DE \}.
\end{eqnarray*}
The vertex classes $\U$ and $\Z$ are mutually exclusive, i.e. $|\U|\cdot|\Z| = 0$, and the matrix $M$ is given by
\[
\bordermatrix{
    & \C & \Cp & \Cm & \I & \U \cr
\C  &  1 &  *  &  1  &  * &  1 \cr
\Cp &  1 &  *  &  1  &  * &  1 \cr
\Cm &  * &  0  &  *  &  0 &  0 \cr
\I  &  * &  0  &  *  &  0 &  0 \cr
\U  &  1 &  0  &  1  &  0 &  0
}
\qquad \mbox{or} \qquad
\bordermatrix{
    & \C & \Cp & \Cm & \I & \Z \cr
\C  &  1 &  *  &  1  &  * &  1 \cr
\Cp &  1 &  *  &  1  &  * &  1 \cr
\Cm &  * &  0  &  *  &  0 &  0 \cr
\I  &  * &  0  &  *  &  0 &  0 \cr
\Z  &  1 &  0  &  1  &  0 &  1
}.
\]
Denote the set of all digraphs $\DG$ that satisfy this construction by $\gT$.
\label{def:gT}
\end{definition}
\begin{figure}[t]
\centering
\includegraphics[width=5in]{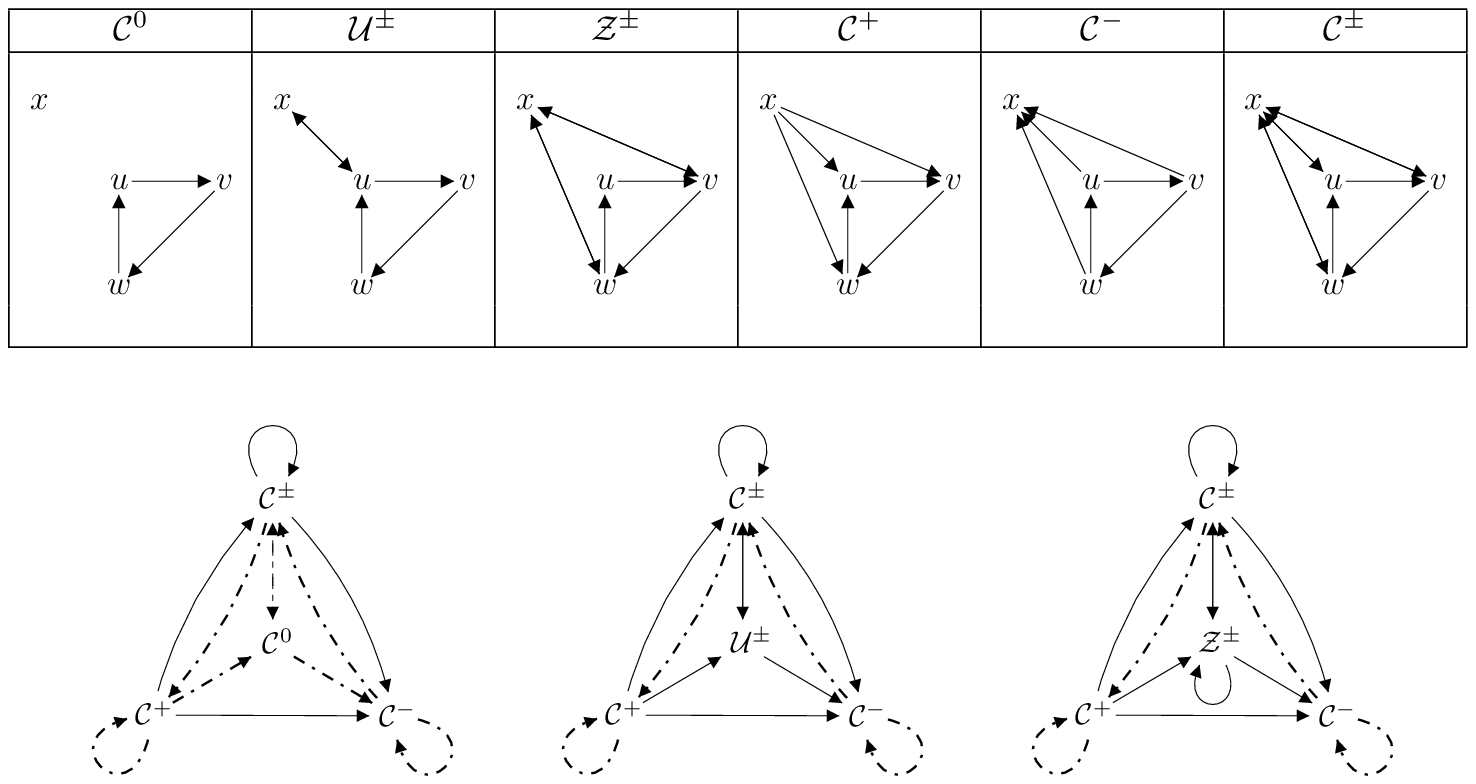}
\caption{{\bf Top}: The six vertex classes $\{\I,\U,\Z,\Cp,\Cm,\C\}$ of $\gT$ defined by how a vertex $x$ in each class connects to a directed 3-cycle $\CC$ denoted by $(u,v,w,u)$.  {\bf Bottom}: Diagrams showing the relations within and between the six possible vertex classes.  Solid and dashed-dotted arrows denote forced and allowable {\de}s, respectively, while the absence of an arrow denotes no {\de}s.  $\C$ and $\Z$ are cliques, $\I$ and $\U$ are independent sets and $G[\Cp]$ and $G[\Cm]$ are arbitrary subgraphs.  $\U$ and $\Z$ are mutually exclusive (i.e. $|\U|\cdot|\Z| = 0$) with no {\de}s between $\I$ and $\U$ or $\I$ and $\Z$.}
\label{fig:3rules}
\end{figure}

Figure~\ref{fig:3rules} shows diagrams depicting the vertex classes and the connections within and between them as determined by the $M$ matrix.  It can be seen almost immediately that if $\DG \in \gT$, then $\DG^{C} \in \gT$ with the following vertex class swaps: $\C \leftrightarrow \I$, $\Cp \leftrightarrow \Cm$, and $\Z \leftrightarrow \U$.  Thus, we move from one $M$-matrix given in Definition~\ref{def:gT} to the other upon complementation.

We now define the set of digraphs $\gQ \subseteq \gR$ as follows:

\begin{definition}
\label{def:gQ}
\[
\gQ \equiv \bigcup_{d\in \gD} \{\DG\in R(d) : \DG \ \mbox{satisfies Lemma~\ref{lem:gQ}}\}.
\]
\end{definition}

\begin{lem}
\label{lem:gQ}
Let $d \in \gD$ with ill-defined maximal index pair $[i,\sK]$ and corresponding vertex pair $[u,\sM] = [\IV{i},\IV{\sK}]$.  There exists a $\DG \in R(d)$ such that the following properties hold:
\begin{tabbing}
\qquad \= (vii) \=\= (a) \= \+ \kill
          (i)   \> $(\Mp,u,\Mm\setminus\Mp) \subseteq P(\DG)$. \\
          (ii)  \> There exists a vertex set $C = \{u,v,w\} \subseteq V(\DG)$ such that $\DG[C] \cong \CC$ \\
                \> with $(u,v,w,u) \in P(\DG)$, $w \in \Mp\cap\Mm$, and $v \notin \Mp\cup\Mm$. \\
          (iii) \> For all $x \in V\setminus C$, \\
                \>       \> (a) $(x,w)\in \DE \Longleftrightarrow (x,v)\in \DE$, \\
                \>       \> (b) $(w,x)\in \DE \Longleftrightarrow (v,x)\in \DE$. \\
          (iv)  \> If $x \in \Mp$ ($\Mm$), then $d_{x} \pgeq d_{w}$ ($d_{x} \ngeq d_{w}$). \\
          (v)   \> If $x,y \in V\setminus C$ and $(x,y)\in A$, then either $(x,v)\in A$ or $(v,y)\in A$. \\
          (vi)  \> For all $x \in V\setminus C$, $(u,x)\in A \Longleftrightarrow x\in \Mm$. \\
          (vii) \> For all $x \in V\setminus C$, \\
                \>       \> (a) if $d_{w}^{-}-1 \leq d_{x}^{-} \leq d_{w}^{-}+1$, then for $y \in V\setminus C$, $y \neq x$, \\
                \>       \>     \> $(y,w)\in \DE \Longleftrightarrow (y,x)\in \DE$, \\
                \>       \> (b) if $d_{w}^{+}-1 \leq d_{x}^{+} \leq d_{w}^{+}+1$, then for $y \in V\setminus C$, $y \neq x$, \\
                \>       \>     \> $(w,y)\in \DE \Longleftrightarrow (x,y)\in \DE$. \\
\end{tabbing}
\end{lem}


\bnprf
The results for $(i)$-$(iii)$ follow from Lemmas~\ref{lem:arbcon} and \ref{lem:vweqz}.
\begin{itemize}
\item[($iv$)]
As was shown in Theorem~\ref{thm:digraphic}, properties $(ii)$ and $(iii)$ imply $d_{v} = d_{w}$, and since $v \notin \Mp\cup\Mm$ and $w \in \Mp\cap\Mm$, we have the result.
\item[($v$)]
Suppose, to the contrary, that $(x,y)\in \DE$ but $\{(x,v),(v,y)\} \subseteq \DE^{C}$.  Then we can perform the 3-switch $\sigma_{3}\bigl((u,v),(v,w),(x,y)\bigr)$.

\item[($vi$)]
Suppose $(u,x)\in \DE$ and $x \notin \Mm$.  By $(v)$, either $(\{v,w\},x)\subseteq \DE^{C}$ or $(\{v,w\},x)\subseteq \DE$.  If $(\{v,w\},x)\subseteq \DE^{C}$, then since $x \notin \Mm$, we can perform the 2-switch $\sigma_{2}\bigl((u,x),(v,w)\bigr)$.  So we must have $(\{v,w\},x)\subseteq \DE$.  Since $x \notin \Mm$, $3 \leq d_{x}^{-} \leq d_{w}^{-}$ and thus there is a $y \notin C \cup \{x\}$ such that $(y,w)\in \DE$ and $(y,x)\notin \DE$.  But then since $x \notin \Mm$, we can perform the 2-switch $\sigma_{2}\bigl((u,x),(y,w)\bigr)$.  Thus, we must have $x \in \Mm$.  For the converse, suppose $x \in \Mm$ and $(u,x)\notin \DE$.  But this means there exists a $z \notin \Mm\cup\{v\}$ such that $(u,z)\in\DE$, which we just showed was a contradiction.

\item[($vii$)]
The proofs of $(a)$ and $(b)$ are analogous, so we will only prove $(a)$.  Suppose first that $(y,w)\in \DE$ and $(y,x)\notin \DE$.  We have $2 \leq d_{w}^{-} \leq d_{x}^{-}+1$, and thus $d_{x}^{-} \geq 1$.  If $(\{v,w\},x)\subseteq \DE$, we have a contradiction by the 3-switch $\sigma_{3}\bigl((u,v),(y,w),(w,x)\bigr)$.  Thus, $(\{v,w\},x)\subseteq \DE^{C}$.  If $d_{x}^{-} = d_{w}^{-}-1 < d_{w}^{-}$ and $(u,x)\in \DE$, then $x \in \Mm$ by $(vi)$.  But by $(iv)$ we have $d_{x}^{-} \geq d_{w}^{-}$, which is a contradiction.  We are left with either $d_{x}^{-} = d_{w}^{-}-1$ and $(u,x)\notin \DE$, or $d_{x}^{-} \geq d_{w}^{-}$.  In either case, we must have a vertex $z \notin C$ such that $(z,x)\in \DE$ and $(z,\{v,w\})\subseteq \DE^{C}$.  But by $(v)$, either $(z,v)\in \DE$ or $(v,x)\in \DE$, both of which can't happen.  Thus, we have a contradiction. \\

Now suppose $(y,x)\in \DE$ and $(y,w)\notin \DE$.  By $(v)$, $(\{v,w\},x)\subseteq \DE \Rightarrow d_{x}^{-} \geq 3$.  If $d_{x}^{-} \leq d_{w}^{-}$, then $d_{w}^{-} \geq 3$ so there is a $z \notin C\cup\{x\}$ such that $(z,w)\in \DE$ and $(z,x)\notin \DE$.  This gives us a contradiction by the 3-switch $\sigma_{3}\bigl((u,v),(z,w),(w,x)\bigr)$.  If $d_{x}^{-} = d_{w}^{-} + 1$, then $d_{x}^{-} > d_{w}^{-}$ and $x \in \Mm$.  By $(vi)$, $(u,x)\in \DE$ and thus $d_{x}^{-} \geq 4 \Rightarrow d_{w}^{-} \geq 3$.  Thus, there is a $z \notin C\cup\{x\}$ such that $(z,w)\in \DE$ and $(z,x)\notin \DE$, which is a contradiction by the same 3-switch as above.
\end{itemize}
\edprf

Let $\DG = (V,\DE) \in \gQ$ with an ill-defined vertex pair $[u,\sM]$.  By Lemma~\ref{lem:gQ}$(ii)$, there is a vertex set $C=\{u,v,w\}$ such that $\DG[C]\cong \CC$.  Let $x \in V\setminus C$ with $\DG_{x}\equiv\DG[V\setminus\{x\}]$ an induced subgraph of $\DG$.  If we can show $\DG_{x} \in \gQ$, then we can use induction to prove $\gQ \subseteq \gT$.  We will do this by showing that $[u,\sM^{\prime}]$ is an ill-defined vertex pair for $\DG_{x}$, where
\[\sM^{\prime} = \{\Mp\setminus\{x\}, \Mm\setminus\{x\}\}.\]

\begin{lem}
If $\DG = (V,\DE) \in\gQ$ and $x\in V\setminus C$, then $\DG_{x}=\DG[V\setminus\{x\}] \in \gQ$.
\label{lem:killx}
\end{lem}
\bnprf

Let $\{\Mp,\Mm\}$ denote the maximal sets for $\DG$ such that there are no arc switches so that $(\Mp,u,\Mm)\subseteq \DP$ and let $w \in \Mp\cap \Mm$.  We will show $\Mp\setminus \{x\}$ and $\Mm \setminus \{x\}$ are valid maximal sets.  Assuming for now that this is true, and supposing $\DG_{x} \notin \gR$, there exists a graph $\DG^{\prime}_{x}\in R(d_{\DG_{x}})$ such that
\[(\Mp\setminus \{x\},u,\Mm\setminus \{x\})\subseteq P(\DG^{\prime}_{x}).\]
Let $\DG^{\prime}$ be a graph such that $V(\DG^{\prime}) = V$ and $\DE(\DG^{\prime}) = \DE(\DG^{\prime}_{x})\cup [\DE(\DG)\setminus \DE(\DG_{x})] \in R(d_{\DG})$.  Thus, we have
\[(\Mp\setminus \{x\},u,\Mm\setminus \{x\})\subseteq \DP(\DG^{\prime}).\]
If $x \in \Mp$, by Lemma~\ref{lem:gQ}$(i)$ and since $\DG\in\gQ$, $(x,u)\in \DE(\DG)\setminus \DE(\DG_{x}) \subseteq \DE(\DG^{\prime})$ (for $x\in\Mm$, $(u,x) \in \DE(\DG^{\prime})$ follows from Lemma~\ref{lem:gQ}$(vi)$).  Thus, we see that we actually have
\[(\Mp,u,\Mm)\subseteq P(\DG^{\prime})\]
which is a contradiction since $\DG \in \gR$.  Thus, $\DG_{x} \in \gR$.  With $x\in V\setminus C$, we see that conditions $(i)$--$(vii)$ in Lemma~\ref{lem:gQ} are satisfied for $\DG_{x}$, showing $\DG_{x} \in \gQ$.

We now show $\Mm\setminus \{x\}$ is a valid maximal in-degree set.  The proof for $\Mp\setminus \{x\}$ is analogous.  Define the following vertex-sets:
\[
\begin{array}{rcl}
W_{-1} & \equiv & \{y \in V\setminus C \ : \ d_{y}^{-} = d_{w}^{-}-1\} \\
W      & \equiv & \{y \in V\setminus C \ : \ d_{y}^{-} = d_{w}^{-}\} \\
W_{+1} & \equiv & \{y \in V\setminus C \ : \ d_{y}^{-} = d_{w}^{-}+1\}.
\end{array}
\]
\begin{figure}[t]
\centering
\includegraphics[width=3in]{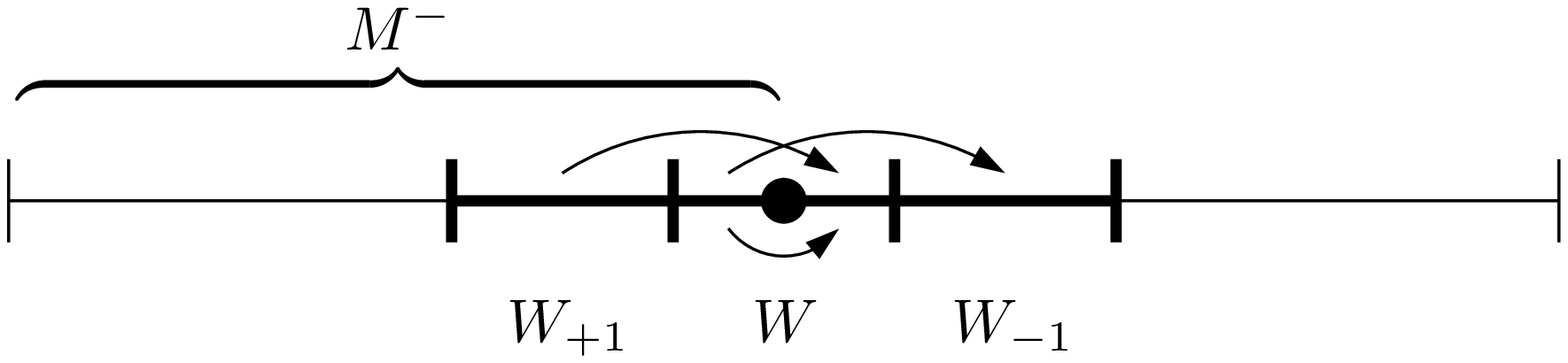}
\caption{Representation of the vertices ordered lexicographically in the negative ordering.  $\Mm$ is a maximal vertex set consisting of vertices with largest degree relative to the negative ordering.  $W_{+1}$, $W$, and $W_{-1}$ are the vertices $y\in V\setminus C$ such that $d_{y}^{-} = d_{w}^{-}+1$, $d_{y}^{-} = d_{w}^{-}$, and $d_{y}^{-} = d_{w}^{-}-1$, respectively.  The bullet denotes the vertex $w$, while the arrows illustrate all possible routes out of $\Mm$ by a decrement of 1 in $d_{y}^{-}$ (top two arrows) or $d_{y}^{+}$ (bottom arrow).}
\label{fig:linegraph}
\end{figure}
Figure \ref{fig:linegraph} displays a line depicting the negative ordering of the vertices.  Note that $w \in W$, which is shown as a $\bullet$ on the line, and is therefore the boundary of $\Mm$ by Lemma~\ref{lem:gQ}$(iv)$.  The arrows denote the only possible fluxes out of $\Mm$ when an {\de } is removed (i.e. an out-degree or in-degree gets decreased by 1).  The top two arrows denote an in-degree decreasing by 1 (i.e. $(x,y)\in \DE, \ y \in W\cup W_{\pm 1}$), while the bottom arrow denotes an out-degree decreasing by 1 ($(y,x)\in \DE$).  But by Lemma~\ref{lem:gQ}$(vii)$, either $(x,W\cup W_{\pm 1})\subseteq \DE$ or $(x,W\cup W_{\pm 1})\subseteq \DE^{C}$.  This shows the in-degrees of all vertices in $W\cup W_{\pm 1}$ decrease by the same amount, and thus no flux out of $\Mm$ via the top arrows is possible.  Similarly, if $y \in W$, the only remaining way to leave $\Mm$ is for $d_{y}^{+} = d_{w}^{+}$ and $(y,x)\in \DE$.  But again Lemma~\ref{lem:gQ}$(vii)$ tells us $(W,x)\subseteq \DE$ and thus there can be no flux out of $\Mm$ via the bottom arrow.  Therefore, we have $\Mm\setminus \{x\}$ is a valid maximal negatively ordered set.
\edprf

We can now prove $\gQ \subseteq \gT$.  First we define
\begin{eqnarray*}
\gQ_{n} & \equiv & \{\DG = (V,\DE)\in \gQ \ : \ |V \setminus C| = n\} \\
\gT_{n} & \equiv & \{\DG = (V,\DE)\in \gT \ : \ |V \setminus C| = n\}
\end{eqnarray*}
such that $\gQ = \cup_{n=0}^{\infty}\gQ_{n}$ and $\gT = \cup_{n=0}^{\infty}\gT_{n}$.
We need a lemma which will be the seed for induction.

\begin{lem}
$\gQ_{i} \subseteq \gT_{i}$, for $i=0,1,2$.
\label{lem:seed}
\end{lem}
\bnprf
[$i=0$]: $\gQ_{0} = \{\CC\} \subseteq \gT_{0}$.

[$i=1$]: Suppose $\DG=(V,\DE)\in\gQ_{1}$ and $x\in V\setminus C$.  We have two types of relations with $C$:  the relation between $x$ with $Z = \{v,w\}$ and $x$ with $u$.  For each type, there are four possibilities.  We introduce the following short-lived notation for its utility in showing in a table all resulting possibilities.
\begin{eqnarray*}
Z^{0} \ & \equiv & \ \{x \in V\setminus C \ : \ (x,Z)\cup (Z,x)\subseteq \DE^{C}\} \\
Z^{+} \ & \equiv & \ \{x \in V\setminus C \ : \ (x,Z)\subseteq \DE \ \mbox{and} \ (Z,x)\subseteq \DE^{C}\} \\
Z^{-} \ & \equiv & \ \{x \in V\setminus C \ : \ (x,Z)\subseteq \DE^{C} \ \mbox{and} \ (Z,x)\subseteq \DE\} \\
Z^{\pm} \ & \equiv & \ \{x \in V\setminus C \ : \ (x,Z)\cup(Z,x)\subseteq \DE\}
\end{eqnarray*}
We define similarly $U^{0}$, $U^{+}$, $U^{-}$ and $U^{\pm}$ for the four possible relations between $x$ and $u$.  The following table displays all possible intersections of these vertex classes.
\[
\begin{array}{c|cccc|}
        & Z^{0}     & Z^{+}     & Z^{-}     & Z^{\pm}   \\ \hline
 U^{0}  & \I        & \cdot & \cdot & \Z        \\
 U^{+}  & \cdot & \Cp       & \cdot & \cdot \\
 U^{-}  & \cdot & \cdot & \Cm       & \cdot \\
U^{\pm} & \U        & \cdot & \cdot & \C        \\ \hline
\end{array}
\]
The vertex classes denoted in the table that are represented in $\gQ_{1}$ are the ones defined in Definition~\ref{def:gT}, while the presence of a $\cdot$ denotes a vertex class that is not represented in $\gQ_{1}$.  Proofs for all cases listed are located in the Appendix in Table~\ref{tab:R1}, giving $\gQ_{1} \subseteq \gT_{1}$.

\begin{figure}[t]
\centering
\includegraphics[width=2.5in]{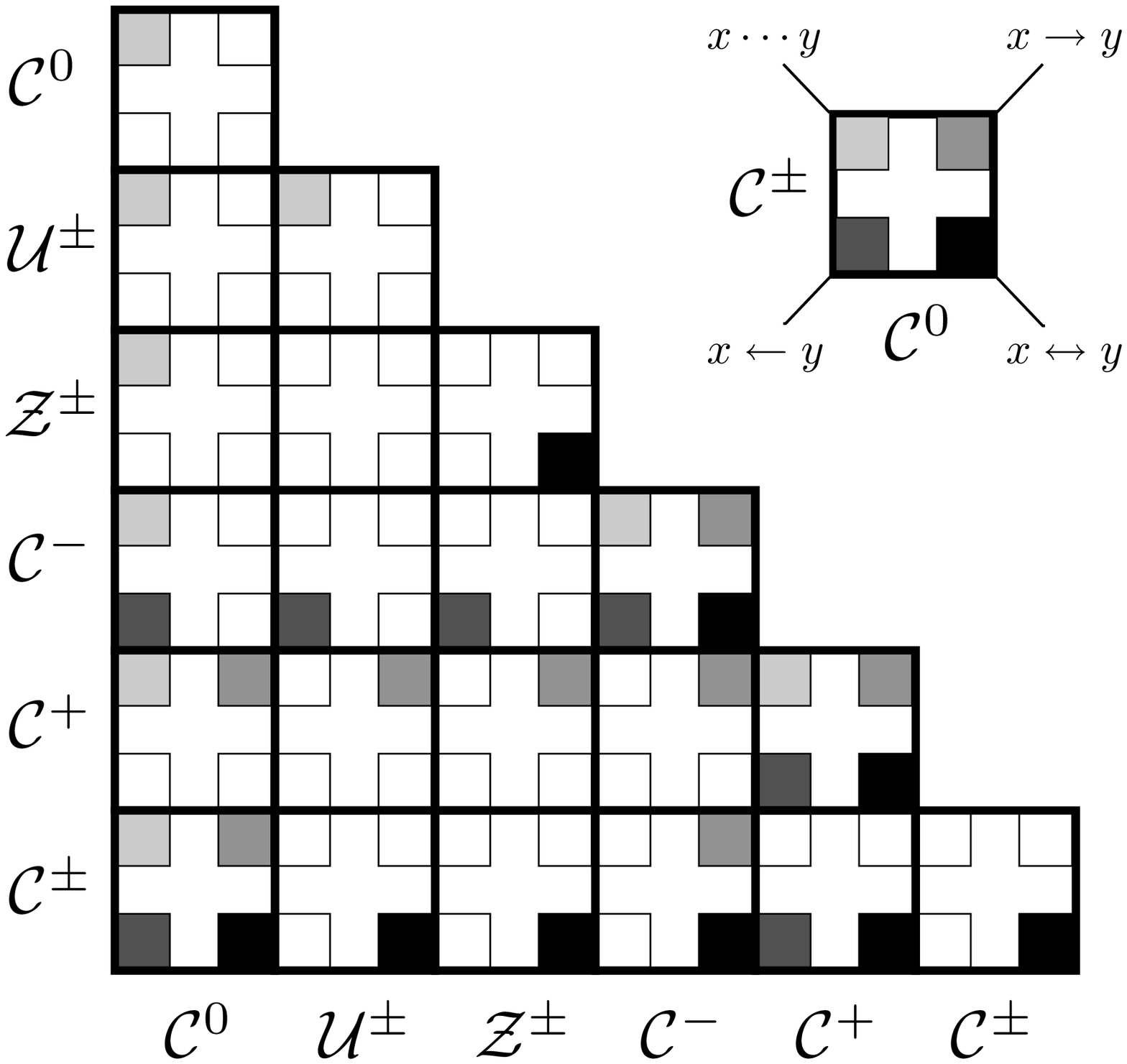}
\caption{Diagram showing the relations within and between vertex classes.  The ``key'' in the upper right shows the $[\C,\I]$ block.  To illustrate, if $x \in \C$ and $y \in \I$, then since the upper right corner of the $[\C,\I]$ block is shaded, $x \rightarrow y$ is an allowable {\de }.  The other corners are detailed in the diagram, with each relation occupying a respective corner and having a specific shade of grey if allowed.  See Table~\ref{tab:R2} for a detailed list of all cases.}
\label{fig:gTable}
\end{figure}
[$i=2$]: Suppose $\DG=(V,\DE)\in \gQ_{2}$ and let $x, y \in V \setminus C$.  By Lemma~\ref{lem:killx}, $\DG_{x} = \DG[V \setminus \{x\}] \in \gQ_{1}$, so $y$ belongs to one of the six classes above.  The same holds for $x$ when $y$ is removed.  Figure~\ref{fig:gTable} displays all possible ways that vertices in the six vertex classes can relate to each other, giving a total of 84 possible relations.  To illustrate, suppose $x \in \C$ and $y \in \Cp$.  The allowable relations are $x \leftarrow y$ (lower-left) or $y \leftrightarrow x$ (lower-right), and thus we can't have $x \rightarrow y$ (upper-right) or no {\de}s between $x$ and $y$ (upper-left).  Thus, if a box is not shaded, then that is not an allowable relation.  Figure~\ref{fig:gTable} corresponds precisely with the $M$-matrix given in Definition~\ref{def:gT}.  Table~\ref{tab:R2} in the Appendix gives proofs for all the cases, showing $\gQ_{2} \subseteq \gT_{2}$.
\edprf

\begin{lem}
$\gQ \subseteq \gT$
\label{lem:gQgT}
\end{lem}
\bnprf
Suppose $\gQ_{k} \subseteq \gT_{k}$.  Let $\DG = (V,\DE) \in \gQ_{k+1}$ and $V\setminus C = \{x_{1},\ldots,x_{k+1}\}$.  We can remove $x_{k+1}$ with the induced subgraph $\DG_{x_{k+1}} = \DG[V\setminus \{x_{k+1}\}]$.  By Lemma~\ref{lem:killx}, $\DG_{x_{k+1}} \in \gQ_{k}$, and so we have $\{x_{1},\ldots,x_{k}\}$ are in the six vertex classes with all {\de}s given by Figure~\ref{fig:gTable}.  If we remove $x_{k}$, we get $x_{k+1}$ in one of the six vertex classes, showing all vertices can be classified {\it and} all {\de}s between $x_{k+1}$ and $x_{i}$, $i=1,\ldots,x_{k-1}$, are given by Figure~\ref{fig:gTable}.  By removing one more vertex, for example $x_{k-1}$, we have the final {\de}s between $x_{k+1}$ and $x_{k}$ given by Figure~\ref{fig:gTable}, and thus $\gQ_{k+1} \subseteq \gT_{k+1}$.  By induction, we have $\gQ \subseteq \gT$.
\edprf


\subsection{Degree sequence characterization}\label{sec:gSgR}

In this section we show that a special set of degree sequences denoted by $\gS$ is a subset of the ill-defined degree sequences $\gD$.

\begin{definition}
Let $d = \{(d^{+}_{i},d^{-}_{i})\}_{i=1}^{N}$ be a degree sequence of length $N$, $J=\{j_{1}, \ldots, j_{n}\}$ a set of indices for $d$ with $3 \leq n \leq N-1$, and $(k,l) \geq (0,0)$ an index pair.  Suppose we have one of the following three cases
\begin{equation}
\begin{array}{rl}
\mbox{\it (i)} & n=3 \ \mbox{and} \ d_{j_1} = d_{j_2} = d_{j_3} = (l+1,k+1) \\
\mbox{\it (ii)} & n > 3 \ \mbox{and} \ d_{j_1} = \cdots = d_{j_{n-1}} = (l+n-2,k+n-2), \ d_{j_{n}} = (l+1,k+1) \\
\mbox{\it (iii)} & n > 3 \ \mbox{and} \ d_{j_1} = (l+n-2,k+n-2), \ d_{j_2} = \cdots = d_{j_{n}} = (l+1,k+1)
\end{array}
\label{equ:gS1}
\end{equation}
with
\begin{equation}
(d_{j_{1}}, \dotsc, d_{j_{n}}) = (\bar{d}_{k+1}, \dotsc, \bar{d}_{k+n}) = (\ubar{d}_{l+1}, \dotsc, \ubar{d}_{l+n})
\label{equ:gS2}
\end{equation}
and the slack sequences satisfying
\begin{equation}
(0, 1, \dotsc, 1, 0) = (\bar{s}_{k}, \bar{s}_{k+1}, \dotsc, \bar{s}_{k+n-1}, \bar{s}_{k+n}) = (\ubar{s}_{l}, \ubar{s}_{l+1}, \dotsc, \ubar{s}_{l+n-1}, \ubar{s}_{l+n}).
\label{equ:gS3}
\end{equation}
Denote the set of all such degree sequences by $\gS$.
\label{def:gS}
\end{definition}

The next theorem proves $\gS \subseteq \gD$ by showing that each case has an ill-defined index for $\HH^{\pm}$.  In particular, we prove that the slack sequence for the residual degree sequence has a negative entry.  
\begin{lem}
$\gS \subseteq \gD$
\label{lem:gSgR}
\end{lem}
\bnprf
{\it (i)--(ii)}  We will show $j_{n}$ is an ill-defined index.  Choose maximal index sets $\sK = \{\Kp,\Km\}$ for $j_{n}$.  Define $\mu$ such that
\[\mu = \sum_{i\in\Kp}\xi^{ij_{n}} + \sum_{i\in\Km}\xi^{j_{n}i}\]
and let $p$ and $q$ denote permutations such that $\bar{d}_{i} = d_{p_{i}}$ and $\ubar{d}_{i} = d_{q_{i}}$.  Let $\nu = \{(\mu_{p_{i}}^{+},\mu_{q_{i}}^{-})\}_{i=1}^{N}$. We have
{\small \setlength\arraycolsep{3pt}
\[
\left(
\begin{array}{l}
\bar{d}^{+} \\
\nu^{+}
\end{array}
\right) = \left(
\begin{array}{ccccccccccc}
\bar{d}^{+}_{1} & \dotsm & \bar{d}^{+}_{k} & l+n-2 & \dotsm  & l+n-2  & l+n-2 & l+1  & \bar{d}^{+}_{k+n+1} & \dotsm & \bar{d}^{+}_{N} \\
-1              & \dotsm & -1                & 0     & \dotsm  & 0      & -1    & -(l+1) & 0                   & \dotsm & 0
\end{array}
\right)
\]
}
as well as
{\small \setlength\arraycolsep{3pt}
\[
\left(
\begin{array}{l}
\ubar{d}^{-} \\
\nu^{-}
\end{array}
\right) = \left(
\begin{array}{cccccccccccc}
\ubar{d}^{-}_{1} & \cdots & \ubar{d}^{-}_{l} & k+n-2 & \cdots  & k+n-2  & k+n-2 & k+1  & \ubar{d}^{-}_{l+n+1} & \cdots & \ubar{d}^{-}_{N} \\
-1               & \cdots & -1                 & 0     & \cdots  & 0      & -1    & -(k+1) & 0                    & \cdots & 0
\end{array}
\right).
\]
}
There are multiple indices where $d_{i} = (l+n-2,k+n-2)$, so we have chosen without loss of generality indices $k+n-1$ and $l+n-1$ of $\nu$ such that $p_{k+n-1} \in \Kp$ and $q_{l+n-1} \in \Km$.  Using the Havel-Hakimi operator $\HH^{\pm}$, define $\tilde{d}$ by
\[\tilde{d} = \Bar{\HH^{\pm}(d,j_{n},\sK)} = \Bar{d-\mu}.\]
The new sequence $\tilde{d}$ is the end result of one simultaneous Havel-Hakimi step, reordered in the positive ordering.  We can also define the new slack sequence by
\[
\tilde{s}_{m} = \sum_{i=1}^{m}[\tilde{d}^{-}]_{i}^{\prime\prime} - \sum_{i=1}^{m}\tilde{d}^{+}_{i}, \quad m=1,\ldots,N.
\]
For the new sequence $\tilde{d}$ to be digraphic, we must have that $\tilde{s}_{m} \geq 0$ for all $m$.  However, we will show below that $\tilde{s}_{k+n-2} = -1$.

Since the slack sequence satisfies $(\ubar{s}_{l}, \ubar{s}_{l+1}, \dotsc, \ubar{s}_{l+n-1}, \ubar{s}_{l+n}) = (0, 1, \dotsc, 1, 0)$ and $[\ubar{d}^{+}]_{l+n}^{\prime\prime}-\ubar{d}^{-}_{l+n} = \ubar{s}_{l+n}-\ubar{s}_{l+n-1}=-1$, we have $[\ubar{d}^{+}]_{l+n}^{\prime\prime} = \ubar{d}^{-}_{l+n}-1 = d^{-}_{j_{n}}-1 = k$.  By the definition of the corrected conjugate sequence, there are $k$ elements in $d^{+}$ such that $d^{+}_{i} \geq l+n-1$.  In particular, we have $\bar{d}^{+}_{i} \geq l+n-1$ for $1 \leq i \leq k$.  Thus, the indices $\{1,\ldots,k+n-2\}$ are only permuted from $\bar{d}$ to $\tilde{d}$, giving
\[
\sum_{i=1}^{k+n-2}\tilde{d}^{+}_{i} = \sum_{i=1}^{k+n-2}\left(\bar{d}^{+}_{i}-\nu^{+}_{i}\right) = \sum_{i=1}^{k+n-2}\bar{d}^{+}_{i} - \sum_{i=1}^{k+n-2}\nu^{+}_{i} = \sum_{i=1}^{k+n-2}\bar{d}^{+}_{i} - k.
\]
We can show similarly for $\bar{d}^{-}$ that since $(\bar{s}_{k}, \bar{s}_{k+1}, \dotsc, \bar{s}_{k+n-1}, \bar{s}_{k+n}) = (0, 1, \dotsc, 1, 0)$ and $[\bar{d}^{-}]_{k+n}^{\prime\prime}-\bar{d}^{+}_{k+n} = \bar{s}_{k+n}-\bar{s}_{k+n-1}=-1$, we have $[\bar{d}^{-}]_{k+n}^{\prime\prime} = \bar{d}^{+}_{k+n}-1 = d^{+}_{j_{n}}-1 = l$.  By the definition of the corrected conjugate sequence, there are $l$ elements in $d^{-}$ such that $d^{-}_{i} \geq k+n-1$.  In particular, $\ubar{d}^{-}_{i} \geq k+n-1$ for $1 \leq i \leq l$.  Each of these gets decreased by $1$ by $\nu^{-}$, but since we are only summing $[\tilde{d}^{-}]^{\prime\prime}$ to $k+n-2 < k+n-1$, we do not see these decreases.  Note also that a permutation of indices from $\bar{d}^{-}$ to $\tilde{d}^{-}$ within $1$ to $k+n-2$ does not change the sum.
Thus, the only elements that get removed from the sum $\sum_{i=1}^{k+n-2}[\bar{d}^{-}]^{\prime\prime}$ are given by the stars below.
\medskip
\begin{center}
\includegraphics[height=0.9in]{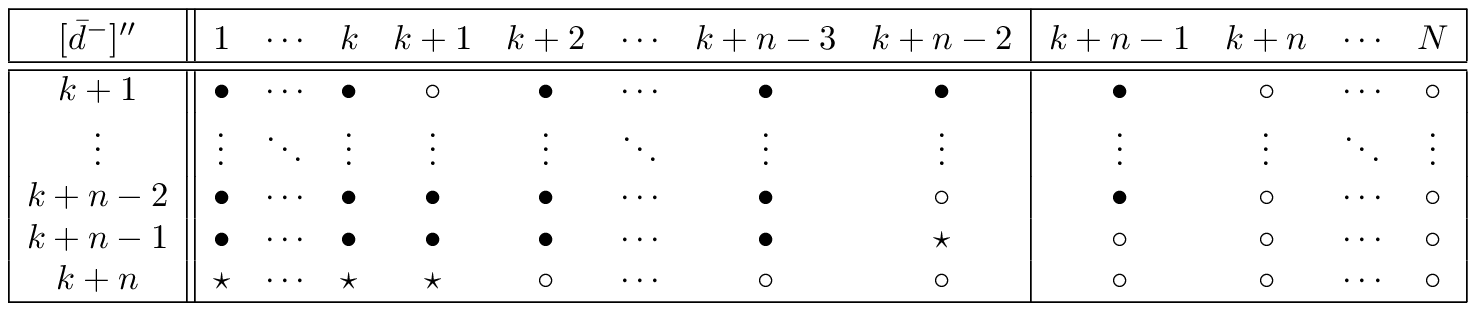}
\end{center}
\medskip
Thus we have
\[
\sum_{i=1}^{k+n-2}[\tilde{d}^{-}]_{i}^{\prime\prime} = \sum_{i=1}^{k+n-2}[\bar{d}^{-}]_{i}^{\prime\prime} - (k+2).
\]
Finally we have
\begin{eqnarray*}
\tilde{s}_{k+n-2} & = & \sum_{i=1}^{k+n-2}[\tilde{d}^{-}]_{i}^{\prime\prime} - \sum_{i=1}^{k+n-2}\tilde{d}^{+}_{i} \\
              & = & \sum_{i=1}^{k+n-2}[\bar{d}^{-}]_{i}^{\prime\prime} - (k+2) - \left(\sum_{i=1}^{k+n-2}\bar{d}^{+}_{i} - k\right) \\
              & = & \left(\sum_{i=1}^{k+n-2}[\bar{d}^{-}]_{i}^{\prime\prime} - \sum_{i=1}^{k+n-2}\bar{d}^{+}_{i}\right) - 2 \\
              & = & \bar{s}_{k+n-2} - 2 \\
              & = & -1
\end{eqnarray*}
since $\bar{s}_{k+n-2} = 1$.

$(iii)$ We will now show $j_{1}$ is an ill-defined index.  Choose maximal index sets $\sK = \{\Kp,\Km\}$ for $j_{1}$.  Define $\mu$ such that
\[\mu = \sum_{i\in\Kp}\xi^{ij_{1}} + \sum_{i\in\Km}\xi^{j_{1}i}\]
and let $p$ and $q$ denote permutations such that $\bar{d}_{i} = d_{p_{i}}$ and $\ubar{d}_{i} = d_{q_{i}}$.  Let $\nu = \{(\mu_{p_{i}}^{+},\mu_{q_{i}}^{-})\}_{i=1}^{N}$. We have $\bar{d}^{+}$ and $\nu^{+}$ similar to above as
{\small \setlength\arraycolsep{3pt}
\[
\left(
\begin{array}{l}
\bar{d}^{+} \\
\nu^{+}
\end{array}
\right) = \left(
\begin{array}{cccccccccccc}
\bar{d}^{+}_{1} & \cdots & \bar{d}^{+}_{k} &   l+n-2  & l+1 &  l+1 & \cdots &  l+1 & \bar{d}^{+}_{k+n+1} & \cdots & \bar{d}^{+}_{N} \\
-1              & \cdots & -1                & -(l+n-2) & 0 & -1 & \cdots & -1 & 0                   & \cdots & 0
\end{array}
\right)
\]
}
as well as $\ubar{d}^{-}$ and $\nu^{-}$ such that
{\small \setlength\arraycolsep{3pt}
\[
\left(
\begin{array}{l}
\ubar{d}^{-} \\
\nu^{-}
\end{array}
\right) = \left(
\begin{array}{cccccccccccc}
\ubar{d}^{-}_{1} & \cdots & \ubar{d}^{-}_{l} &   k+n-2  & k+1 &  k+1 & \cdots &  k+1 & \ubar{d}^{-}_{l+n+1} & \cdots & \ubar{d}^{-}_{N} \\
-1               & \cdots & -1                 & -(k+n-2) & 0 & -1 & \cdots & -1 & 0                    & \cdots & 0
\end{array}
\right).
\]
}
Now there are multiple indices where $d_{i} = (l+1,k+1)$, so we choose without loss of generality indices $\{k+3,\ldots,k+n\}$ and $\{l+3,\ldots,l+n\}$ of $\nu$ such that $\{p_{k+3},\ldots,p_{k+n}\} \subseteq \Kp$ and $\{q_{l+3},\ldots,q_{l+n}\} \subseteq \Km$.  Again define $\tilde{d}$ by
\[\tilde{d} = \Bar{\HH^{\pm}(d,j_{1},\sK)} = \Bar{d-\mu}.\]
We now show that the slack sequence is negative in the $(k+1)$-st element, i.e. $\tilde{s}_{k+1} = -1$.

Since the slack sequence satisfies $(\bar{s}_{k}, \bar{s}_{k+1}, \dotsc, \bar{s}_{k+n-1}, \bar{s}_{k+n}) = (0, 1, \dotsc, 1, 0)$, we have $[\bar{d}^{-}]_{k+n}^{\prime\prime} = \bar{d}^{+}_{k+n}-1 = d^{+}_{j_{n}}-1 = l$.  Thus, there are $l$ elements in $d^{-}$ such that $d^{-}_{i} \geq k+n-1$.  In particular, $\ubar{d}^{-}_{i} \geq k+n-1$ for $1 \leq i \leq l$.  Each of these gets decreased by $1$ by $\nu^{-}$, but since we are only summing $[\tilde{d}^{-}]^{\prime\prime}$ to $k+1 < k+n-1$, we do not see these decreases.  However, unlike before, we also have $\nu^{+}_{i} = -1$ for $i=k+3,\ldots,k+n$ where $\bar{d}^{-}_{i} = k+1$.  We will thus see these $n-2$ decreases.  We will also see a decrease in the $(k+1)$-st position where $\bar{d}^{-}_{k+1} = k+n-2$.  Thus, the elements that get removed from the sum $\sum_{i=1}^{k+1}[\bar{d}^{-}]^{\prime\prime}$ are given by the stars and asterisk below in columns 1 through $k+1$:
\medskip
\begin{center}
\includegraphics[height=0.9in]{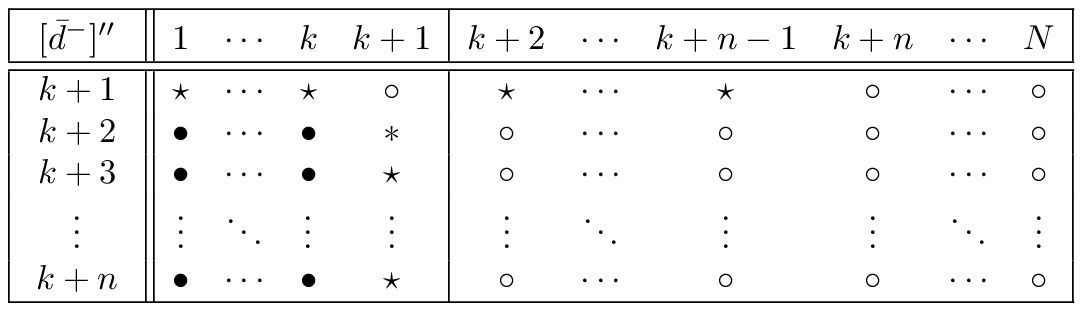}
\end{center}
\medskip
Since we are removing an index within the range of summation $\{1,\ldots,k+1\}$, there is a hidden decrement of $1$ when $\bar{d}^{-}_{k+2}$ moves into the $(k+1)$-st row of $\tilde{d}^{-}$.  This action removes 1 unit since the element that was in the $(k+1)$-st column (see $\ast$ in the above diagram) gets moved to the $(k+2)$-nd column.  Thus, totalling it all up, we have
\[
\sum_{i=1}^{k+1}[\tilde{d}^{-}]_{i}^{\prime\prime} = \sum_{i=1}^{k+1}[\bar{d}^{-}]_{i}^{\prime\prime} - k - (n-2) - 1 = \sum_{i=1}^{k+1}[\bar{d}^{-}]_{i}^{\prime\prime} - (k+n) + 1.
\]
Now, for $\tilde{d}^{+}$ we have
\begin{eqnarray*}
\sum_{i=1}^{k+1}\tilde{d}^{+}_{i} & = & \sum_{i=1}^{k}\tilde{d}^{+}_{i} + \tilde{d}^{+}_{k+1} \\
                                & = & \sum_{i=1}^{k}\bar{d}^{+}_{i} - k + \tilde{d}^{+}_{k+1} \\
                                & = & \sum_{i=1}^{k+1}\bar{d}^{+}_{i} - k + \left[\tilde{d}^{+}_{k+1} - \bar{d}^{+}_{k+1}\right] \\
                                & = & \sum_{i=1}^{k+1}\bar{d}^{+}_{i} - k + \left[l+1 - (l+n-2)\right] \\
                                & = & \sum_{i=1}^{k+1}\bar{d}^{+}_{i} - (k+n) + 3.
\end{eqnarray*}
The second line follows since, as in part $(i)$--$(ii)$, we have $\bar{d}^{+}_{i} \geq l+n-2 > l+2$ for $i=1,\ldots,k$ and thus $\tilde{d}^{+}_{i} \geq l+2$, showing that the ordering is preserved within $1$ to $k$ from $\bar{d}^{+}$ to $\tilde{d}^{+}$.  Subtracting the equalities, we have
\begin{eqnarray*}
\tilde{s}_{k+1} & = & \sum_{i=1}^{k+1}[\tilde{d}^{-}]_{i}^{\prime\prime} - \sum_{i=1}^{k+1}\tilde{d}^{+}_{i} \\
              & = & \sum_{i=1}^{k+1}[\bar{d}^{-}]_{i}^{\prime\prime} - (k+n) + 1 - \left(\sum_{i=1}^{k+1}\bar{d}^{+}_{i} - (k+n) + 3\right) \\
              & = & \left(\sum_{i=1}^{k+1}[\bar{d}^{-}]_{i}^{\prime\prime} - \sum_{i=1}^{k+1}\bar{d}^{+}_{i}\right) - 2 \\
              & = & \bar{s}_{k+1} - 2 \\
              & = & -1
\end{eqnarray*}
since $\bar{s}_{k+1} = 1$.
\edprf


\subsection{Equivalence of characterizations}\label{sec:gTgS}

This section consists of two lemmas that show inclusion for the characterizations in Definitions~\ref{def:gT} and \ref{def:gS}, in particular $d_{\gT} \subseteq \gS$ and $R(\gS) \subseteq \gT$.  These inclusions are used in Theorem~\ref{thm:charequ} to show that they are in fact equivalent, i.e. $d_{\gT} = \gS$ or $\gT = R(\gS)$.

We prove in the next lemma that $d_{\gT} \subseteq \gS$.  Based on the structure of $\gT$, we quickly derive a set of inequalities for the vertex classes that gives us Eqs.~(\ref{equ:gS1}) and (\ref{equ:gS2}) of Definition~\ref{def:gS}.  We then use inherent similarities between the corrected Ferrers diagram with the adjacency matrix to show Eq.~(\ref{equ:gS3}).
\begin{lem}
$d_{\gT} \subseteq \gS$
\label{lem:gTgS}
\end{lem}
\bnprf
Let $\DG = (V,\DE) \in \gT$ and suppose that $\U = \emptyset$.  Define the following constants
\[
\begin{array}{lcllcl}
i^{\prime} & = & |\I|,  & \qquad j^{\prime} & = & |\C|, \\
n^{\prime} & = & |\Z|,  & \qquad n          & = & n^{\prime}+3, \\
k^{\prime} & = & |\Cp|, & \qquad k          & = & k^{\prime}+j^{\prime}, \\
l^{\prime} & = & |\Cm|, & \qquad l          & = & l^{\prime}+j^{\prime}. \\
\end{array}
\]
The constraints on the degrees for each vertex class are listed in Table \ref{tab:ineq}, which are found from Figure~\ref{fig:3rules} and the constants above.  From Table \ref{tab:ineq}, the positive and negative orderings of the vertex classes are given by
\begin{table}[t]
\[
\begin{array}{|ccc||ccc|ccc|}
\hline
& x \in \Z      &&& d_{x}^{+} = l+n-2             &&& d_{x}^{-} = k+n-2             & \\ \hline
& x \in \{v,w\} &&& d_{x}^{+} = l+n-2             &&& d_{x}^{-} = k+n-2             & \\ \hline
& x = u         &&& d_{x}^{+} = l+1               &&& d_{x}^{-} = k+1               & \\ \hline\hline
& x \in \U      &&& d_{x}^{+} = l+1               &&& d_{x}^{-} = k+1               & \\ \hline
& x \in \{v,w\} &&& d_{x}^{+} = l+1               &&& d_{x}^{-} = k+1               & \\ \hline
& x = u         &&& d_{x}^{+} = l+n-2             &&& d_{x}^{-} = k+n-2             & \\ \hline\hline
& x \in \I      &&& 0 \leq d_{x}^{+} \leq l       &&& 0 \leq d_{x}^{-} \leq k       & \\ \hline
& x \in \C      &&& l+n-1 \leq d_{x}^{+} \leq N-1 &&& k+n-1 \leq d_{x}^{-} \leq N-1 & \\ \hline
& x \in \Cp     &&& l+n \leq d_{x}^{+} \leq N-1   &&& 0 \leq d_{x}^{-} \leq k-1     & \\ \hline
& x \in \Cm     &&& 0 \leq d_{x}^{+} \leq l-1     &&& k+n \leq d_{x}^{-} \leq N-1   & \\ \hline
\end{array}
\]
\caption{Degree inequalities for the vertex classes $\{\I,\U,\Z,\Cp,\Cm,\C\}$ following directly from the forced and allowable {\de}s in Figure~\ref{fig:3rules}.}
\label{tab:ineq}
\end{table}
\begin{equation}
\begin{array}{l}
d_{\C \cup \Cp} \pgrt d_{w} = d_{v} = d_{\Z} \pgeq d_{u} \pgrt d_{\Cm \cup \I} \\
d_{\C \cup \Cm} \ngrt d_{w} = d_{v} = d_{\Z} \ngeq d_{u} \ngrt d_{\Cp \cup \I}.
\end{array}
\label{eq:classorder}
\end{equation}
For $\Z = \{z_{1},\ldots,z_{n-3}\}$, let
\[(d_{1}, \ldots, d_{n}) = (d_{v}, d_{w}, d_{z_{1}}, \ldots, d_{z_{n-3}}, d_{u}).\]
Note that $(k,l) \geq (0,0)$ is an index pair with
\[(d_{1}, \ldots, d_{n}) = (\bar{d}_{k+1}, \ldots, \bar{d}_{k+n}) = (\ubar{d}_{l+1}, \ldots, \ubar{d}_{l+n}).\]
We also have, as in case $(i)$ in Definition~\ref{def:gS},
\[d_{1} = \cdots = d_{n-1} = (l+n-2,k+n-2), \ d_{n} = (l+1,k+1).\]
It remains for us to show that
\[(0, 1, \dotsc, 1, 0) = (\bar{s}_{k}, \bar{s}_{k+1}, \dotsc, \bar{s}_{k+n-1}, \bar{s}_{k+n}) = (\ubar{s}_{l}, \ubar{s}_{l+1}, \dotsc, \ubar{s}_{l+n-1}, \ubar{s}_{l+n}).\]
\begin{figure}
\centering
\includegraphics[width=5in]{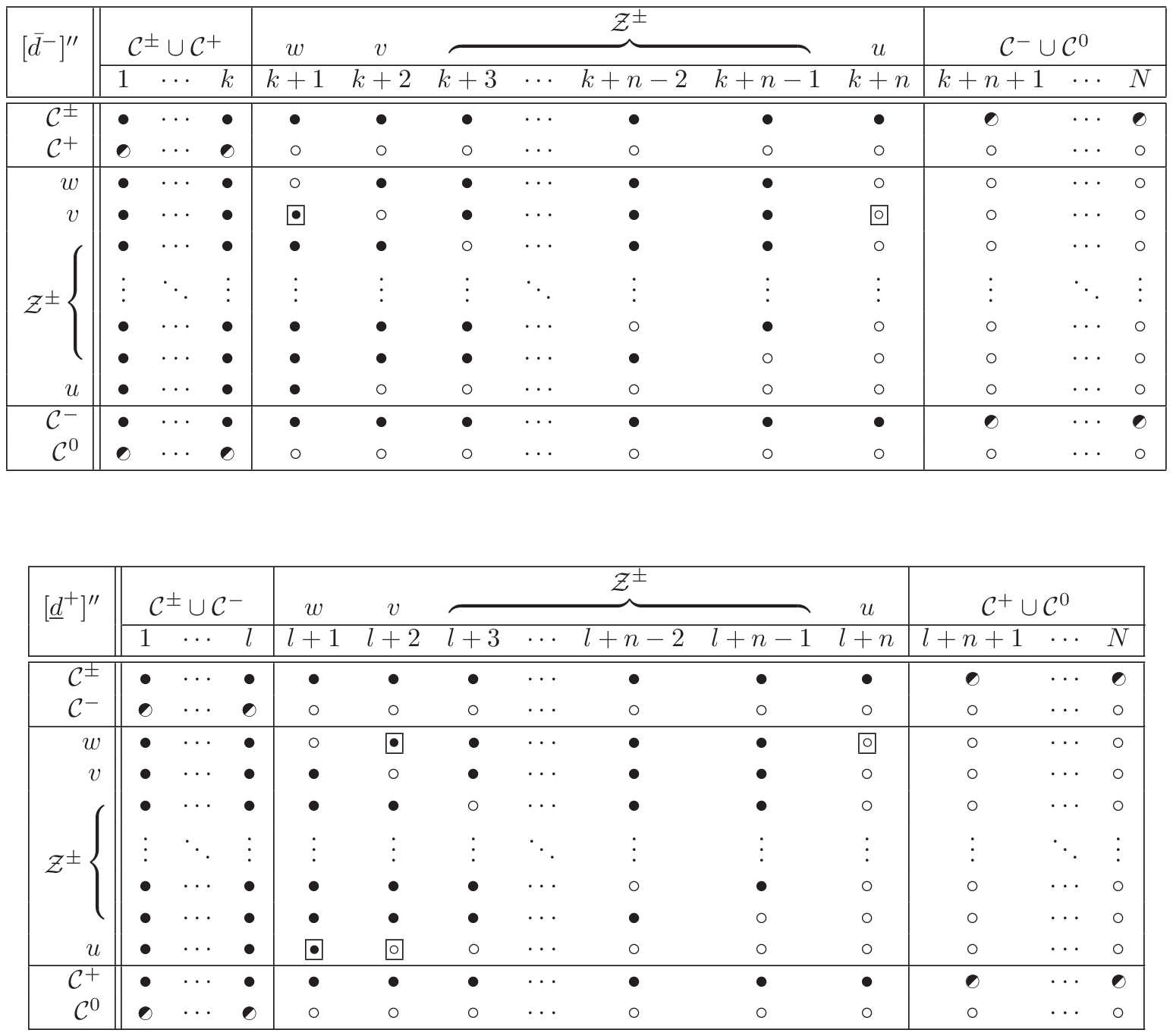}
\caption{The corrected Ferrers diagram for $\bar{d}^{-}$ and $\ubar{d}^{+}$ when $|\U|=\emptyset$.  Half-closed circles denote places where there may or may not be a closed circle.  In the top (bottom) diagram, all closed and open circles correspond to a 1 or 0, respectively, in the transposed (untransposed) adjacency matrix for every digraph $\DG\in R(d)$.  The only locations where there is not a correspondence, besides possibly where half-closed circles are located, are denoted by the boxed-closed and boxed-open circles where there is a 0 and 1, respectively, in the transposed (top diagram) or untransposed (bottom diagram) adjacency matrix.}
\label{fig:ferZ}
\end{figure}
\begin{figure}
\centering
\includegraphics[width=5in]{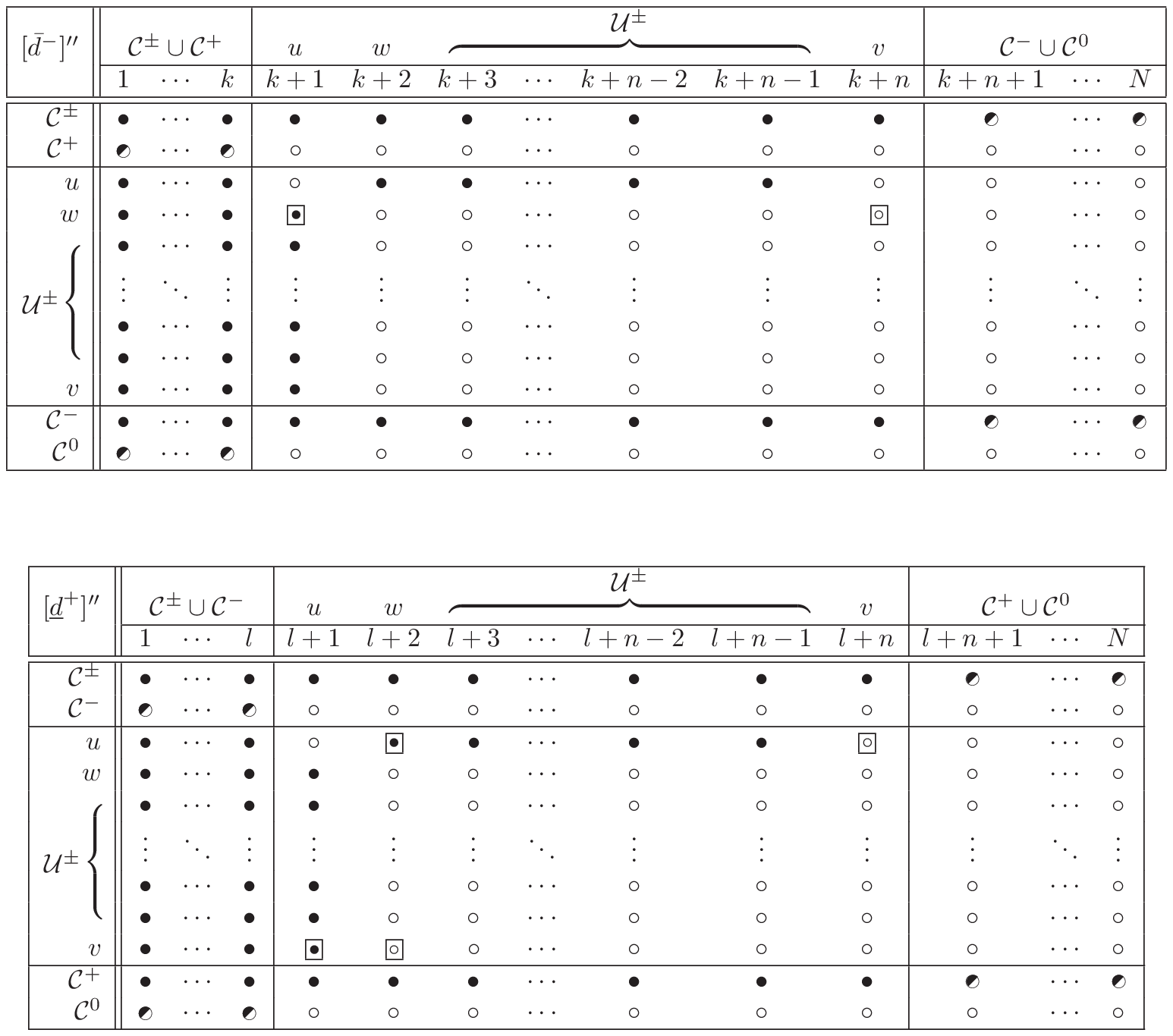}
\caption{The corrected Ferrers diagram for $\bar{d}^{-}$ and $\ubar{d}^{+}$ when $|\Z|=\emptyset$.  Description same as in Figure~\ref{fig:ferZ} caption.}
\label{fig:ferU}
\end{figure}
Given the inequalities in Table \ref{tab:ineq}, we have drawn the corrected Ferrers diagrams for $\bar{d}^{-}$ and $\ubar{d}^{+}$ in Fig.~\ref{fig:ferZ}, where half-closed circles denote locations where there may be a closed or open circle.  Note that the corrected Ferrers diagrams for $\ubar{d}^{+}$ and $\bar{d}^{-}$ correspond closely with the adjacency matrix and its transpose, respectively, for $\DG$, where 1's and 0's are in place of closed and open circles.  We will consider $\bar{d}^{-}$ in regards to $\bar{s}$, with an analogous argument for $\ubar{d}^{+}$ regarding $\ubar{s}$.  There are only two places where there is a definite mismatch between the corrected Ferrers diagram for $\bar{d}^{-}$ and the transposed adjacency matrix, which are labeled with a box surrounding the circles.  The boxes around the closed and open circles correspond in this diagram to the fact that $(w,v)\notin \DE$ and $(u,v)\in \DE$, respectively.  The other locations where there may be a mismatch occur at the half-closed circles.

The corrected Ferrers diagram efficiently illustrates how the slack sequence satisfies the constraints in Definition~\ref{def:gS}.  For the first $k$ columns corresponding to the set $\C\cup\Cp$, the locations where there are mismatches between the corrected Ferrers diagram and the transposed adjacency matrix occur at the rows corresponding to $\Cp$ and $\I$.  However, note from Fig.~\ref{fig:3rules} that the only possible {\de}s into $\Cp$ or $\I$ come from $\C\cup\Cp$, and in the rows of the corrected Ferrers diagram corresponding to $\Cp$ and $\I$, the only places where a filled circle could be is in the $\C\cup\Cp$ columns.  This shows that the number of filled circles in the first $k$ columns of the corrected Ferrers diagram equals the number of arcs out of $\C\cup\Cp$, which means
\[\sum_{i=1}^{k}[\bar{d}^{-}]_{i}^{\prime\prime} = \sum_{i=1}^{k}\bar{d}_{i}^{+},\]
i.e. $\bar{s}_{k} = 0.$  As we move to the $(k+1)$-st index, there is a box surrounding a closed circle implying $(w,v)\notin \DE$, which gives $\bar{s}_{k+1} = 1$.  For the columns corresponding to vertices in $\{v,\Z\}$, the corrected Ferrers diagram and transposed adjacency matrix agree, and so $\bar{s}_{i} = 1$ for $i=k+2,\ldots,k+n-1$.  We finally get a decrement in the partial sums at the $(k+n)$-th index with a box surrounding an open circle implying $(u,v)\in \DE$, giving $\bar{s}_{k+n} = 0$.  A similar proof holds for the case $\Z = \emptyset$ by comparing with case $(ii)$ in Definition~\ref{def:gS} (See Fig.~\ref{fig:ferU} for the corresponding corrected Ferrers diagrams).
\edprf

The final theorem needed to prove the equivalence of the ill-defined sets and the characterizations is to show $R(\gS) \subseteq \gT$.  Based on the $M$-partition structure of $\gT$, we need to be able to identify the vertex classes solely from the degree sequences in $\gS$ and show that the degree sequence structure forces the existence or absence of an arc in their digraph realizations.  To accomplish this, we use the following lemma, which follows immediately from proof by contrapositive.
\begin{lem}
\label{lem:FFarc}
Let $d$ be digraphic, $(i,j)$ an index pair such that $i\neq j$ with $(u,v) = (\IV{i},\IV{j})$, and $\xi^{ij} \equiv \{(\delta_{ik},\delta_{jk})\}_{k=1}^{N}$ with $\delta_{ij}$ the Kronecker delta operator.

(i) If $d-\xi^{ij}$ is not digraphic, then $(u,v) \notin \DE(\DG)$ for all $\DG \in R(d)$.

(ii) If $d+\xi^{ij}$ is not digraphic, then $(u,v) \in \DE(\DG)$ for all $\DG \in R(d)$.
\end{lem}

\begin{lem}
\label{lem:gRgT}
\[R(\gS)\subseteq \gT\]
\end{lem}

\bnprf
Let $d \in \gS$.  Using the slack sequences, we have
\begin{eqnarray*}
[\ubar{d}^{+}]^{\prime\prime}_{l+n} - \ubar{d}^{-}_{l+n} = \ubar{s}_{l+n}-\ubar{s}_{l+n-1} = -1 & \Rightarrow & [\ubar{d}^{+}]^{\prime\prime}_{l+n} = \ubar{d}^{-}_{l+n}-1 = d^{-}_{j_{n}}-1 = k \\
& \Rightarrow & \bar{d}^{+}_{i} \geq l+n-1 \ \mbox{for} \ i=1,\ldots,k.
\end{eqnarray*}
We also have
\begin{eqnarray*}
[\ubar{d}^{+}]^{\prime\prime}_{l+1} - \ubar{d}^{-}_{l+1} = \ubar{s}_{l+1}-\ubar{s}_{l} = 1 & \Rightarrow & [\ubar{d}^{+}]^{\prime\prime}_{l+1} = \ubar{d}^{-}_{l+1}+1 = d^{-}_{j_{1}}+1 = k+n-1 \\
& \Rightarrow & \exists \ \mbox{at most} \ k+n-1 \ \mbox{indices such that} \ d^{+}_{i} \geq l+1 \\
& \Rightarrow & \exists \ N-(k+n) \ \mbox{indices such that} \ d^{+}_{i} \leq l \\ 
& \Rightarrow & \bar{d}^{+}_{i} \leq l \ \mbox{for} \ i=k+n+1,\ldots,N.
\end{eqnarray*}
Similar arguments show $\ubar{d}^{-}_{i} \geq k+n-1$ for $i=1,\ldots,l$ and $\ubar{d}^{-}_{i} \leq k$ for $i=l+n+1,\ldots,N$.  Define the following index sets
\begin{eqnarray*}
X^{\pm} \ & \equiv & \ \{ i \ : \ d^{+}_{i} \geq l+n-1, \ d^{-}_{i} \geq k+n-1 \} \\
X^{0} \ & \equiv & \ \{ i \ : \ d^{+}_{i} \leq l, \ d^{-}_{i} \leq k \} \\
X^{+} \ & \equiv & \ \{ i \ : \ d^{+}_{i} \geq l+n-1, \ d^{-}_{i} \leq k \} \\
X^{-} \ & \equiv & \ \{ i \ : \ d^{+}_{i} \leq l, \ d^{-}_{i} \geq k+n-1 \} \\
S     \ & \equiv & \ \{ j_{1}, \ldots, j_{n} \}
\end{eqnarray*}
with the corresponding vertex sets $\C = \IV{X^{\pm}}$, $\I = \IV{X^{0}}$, $\Cp = \IV{X^{+}}$, $\Cm = \IV{X^{-}}$ and $\mathcal{V} = \IV{S}$.  Since $l+1 \leq d_{S}^{+} \leq l+n-2$ and $k+1 \leq d_{S}^{-} \leq k+n-2$, these index sets constitute a partitioning of $\{1,\ldots,N\}$, with their relative ordering given by
\[
\begin{array}{l}
d_{X^{\pm} \cup X^{+}} \pgrt d_{S} \pgrt d_{X^{-} \cup X^{0}} \\
d_{X^{\pm} \cup X^{-}} \ngrt d_{S} \ngrt d_{X^{+} \cup X^{0}}.
\end{array}
\]

If $i \in X^{\pm}\cup X^{+}$ and $j \in X^{\pm}\cup X^{-} \cup S$, then $d^{+}_{i} \geq l+n-1$ and $d^{-}_{j} \geq k+1$.  Consider $\tilde{d} = \overline{d+\xi^{ij}}$.  By the ordering of the degrees, the indices $\{1,\ldots,k\}$ are permuted from $\bar{d}$ to $\tilde{d}$.  Thus,
\begin{eqnarray*}
\tilde{s}_{k} & = & \sum_{m=1}^{k}[\tilde{d}^{-}]^{\prime\prime}_{m} - \sum_{m=1}^{k}\tilde{d}^{+}_{m} \\
                & = & \sum_{m=1}^{k}[\bar{d}^{-}]^{\prime\prime}_{m} - \left(\sum_{m=1}^{k}\bar{d}^{+}_{m} + 1\right) \\
                & = & \bar{s}_{k} - 1 \\
                & = & -1,
\end{eqnarray*}
i.e. $\tilde{d}$ is not digraphic.  By Lemma~\ref{lem:FFarc}$(ii)$, $(\IV{i},\IV{j})$ is a forced arc.  Similarly, let $i \in S$ and $j \in X^{\pm}\cup X^{-}$, and thus $d^{+}_{i} \geq l+1$ and $d^{-}_{j} \geq k+n-1$.  If we consider $\utilde{d} = \underline{d+\xi^{ij}}$, then we have
\begin{eqnarray*}
\utilde{s}_{l} & = & \sum_{m=1}^{l}[\utilde{d}^{+}]^{\prime\prime}_{m} - \sum_{m=1}^{l}\utilde{d}^{-}_{m} \\
                & = & \sum_{m=1}^{l}[\ubar{d}^{+}]^{\prime\prime}_{m} - \left(\sum_{m=1}^{l}\ubar{d}^{-}_{m} + 1\right) \\
                & = & \ubar{s}_{l} - 1 \\
                & = & -1.
\end{eqnarray*}
By Lemma~\ref{lem:FFarc}$(ii)$, $(\IV{i},\IV{j})$ is a forced arc.

Now let $i \in X^{0}\cup X^{-}$ and $j \in X^{0}\cup X^{+} \cup S$.  Thus, we have $d^{+}_{i} \leq l$ and $d^{-}_{j} \leq k+n-2$.  Consider $\tilde{d} = \overline{d-\xi^{ij}}$.  Again, we have
\begin{eqnarray*}
\tilde{s}_{k+n} & = & \sum_{m=1}^{k+n}[\tilde{d}^{-}]^{\prime\prime}_{m} - \sum_{m=1}^{k+n}\tilde{d}^{+}_{m} \\
                & = & \left(\sum_{m=1}^{k+n}[\bar{d}^{-}]^{\prime\prime}_{m}-1\right) - \sum_{m=1}^{k+n}\bar{d}^{+}_{m} \\
                & = & \bar{s}_{k+n} - 1 \\
                & = & -1,
\end{eqnarray*}
i.e. $\tilde{d}$ is not digraphic.  By Lemma~\ref{lem:FFarc}$(i)$, $(\IV{i},\IV{j})$ is a forbidden arc.  Similarly, let $i \in S$ and $j \in X^{+}\cup X^{0}$, and thus $d^{+}_{i} \leq l+n-2$ and $d^{-}_{j} \leq k$.  If we consider $\utilde{d} = \underline{d-\xi^{ij}}$, then
\begin{eqnarray*}
\utilde{s}_{l+n} & = & \sum_{m=1}^{l+n}[\utilde{d}^{+}]^{\prime\prime}_{m} - \sum_{m=1}^{l+n}\utilde{d}^{-}_{m} \\
                & = & \left(\sum_{m=1}^{l+n}[\ubar{d}^{+}]^{\prime\prime}_{m}-1\right)
                - \sum_{m=1}^{l+n}\ubar{d}^{-}_{m} \\
                & = & \ubar{s}_{l+n} - 1 \\
                & = & -1.
\end{eqnarray*}
By Lemma~\ref{lem:FFarc}$(i)$, $(\IV{i},\IV{j})$ is a forbidden arc.

Summarizing our results so far, we have found index sets corresponding to the vertex sets in Definition~\ref{def:gT} and have shown
\[
\bordermatrix{
             & \C & \Cp & \Cm & \I & \mathcal{V} \cr
\C           &  1 &  *  &  1  &  * &  1 \cr
\Cp          &  1 &  *  &  1  &  * &  1 \cr
\Cm          &  * &  0  &  *  &  0 &  0 \cr
\I           &  * &  0  &  *  &  0 &  0 \cr
\mathcal{V}  &  1 &  0  &  1  &  0 &  ?
}.
\]
To complete the picture, we need to show there is a directed 3-cycle $\CC \subseteq \mathcal{V}$, as well as distinguishing the sets $\U$ and $\Z$.  Using the constants from Lemma~\ref{lem:gTgS}, from the definitions of the index sets we have $k = |\C\cup\Cp|$ and $l = |\C\cup\Cm|$.  Connecting these forced arcs, we have for $j_{m} \in S$, $d^{\prime}_{j_{m}} = d_{j_{m}}-(l,k)$.  Thus, the residual sequence $d^{\prime}_{S}$ becomes
\[
\begin{array}{rl}
\mbox{\it (i)} & n=3 \ \mbox{and} \ d^{\prime}_{j_1} = d^{\prime}_{j_2} = d^{\prime}_{j_3} = (1,1) \\
\mbox{\it (ii)} & n > 3 \ \mbox{and} \ d^{\prime}_{j_1} = \cdots = d^{\prime}_{j_{n-1}} = (n-2,n-2), \ d^{\prime}_{j_{n}} = (1,1) \\
\mbox{\it (iii)} & n > 3 \ \mbox{and} \ d^{\prime}_{j_1} = (n-2,n-2), \ d^{\prime}_{j_2} = \cdots = d^{\prime}_{j_{n}} = (1,1).
\end{array}
\]
In case $(i)$, there is a directed 3-cycle $\CC$ with $\Z = \U = \emptyset$, and thus $R(d) \in \gT$.  Suppose now we are in case $(ii)$ so that we have
\[
d^{\prime}_{S} = \left(\begin{array}{cccc} n-2 & \cdots & n-2 & 1 \\ n-2 & \cdots & n-2 & 1\end{array}\right).
\]
Our choice of maximal index sets $\sK$ for index $j_{n}$ can be any of the remaining indices, and so choose $i,j \in \{j_{1},\ldots,j_{n-1}\}$ with $\Kp = \{i\}$ and $\Km = \{j\}$.  Letting $\hat{d}_{S} = \HH^{\pm}(d^{\prime}_{S},j_{n},\sK)$, if we have $i=j$, then it can be easily seen that a loop will be forced, so we must have $i\neq j$.  This reiterates what we already know, in particular $j_{n}$ is an ill-defined index.  Without loss of generality, let $i=j_{n-2}$ and $j=j_{n-1}$.  We now have
\[
\hat{d}_{S} = \left(\begin{array}{cccccc} n-2 & \cdots & n-2 & n-3 & n-2 & 0 \\ n-2 & \cdots & n-2 & n-2 & n-3 & 0 \end{array}\right).
\]
If we define $\Z = \IV{\{j_{1},\ldots,j_{n-3}\}}$, $w=\IV{j_{n-2}}$, $v=\IV{j_{n-1}}$ and $u=\IV{j_{n}}$, then it can be easily seen from the degree sequence $\hat{d}_{S}$ and the connections already made that we must have a directed 3-cycle $(u,v,w,u)$ with bidirectional arcs between $Z=\{v,w\}$ and $\Z$, with $\Z$ a clique.  Thus, $R(d)\in \gT$.

Now suppose we are in case $(iii)$ so that we have
\[
d^{\prime}_{S} = \left(\begin{array}{cccc} n-2 & 1 & \cdots & 1 \\ n-2 & 1 & \cdots & 1\end{array}\right).
\]
Again our choice of maximal index sets $\sK$ for index $j_{1}$ can be any combination of the remaining indices.  However, it is quickly seen that if we have $\Kp = \Km$, then we will have a loop, showing again that in this case $j_{1}$ is ill-defined.  So, without loss of generality, let $\Kp = \{j_{2},\ldots,j_{n-1}\}$ and $\Km = \{j_{2},\ldots,j_{n-2},j_{n}\}$.  Setting $\hat{d}_{S} = \HH^{\pm}(d^{\prime}_{S},j_{1},\sK)$, we have
\[
\hat{d}_{S} = \left(\begin{array}{ccccc} 0 & \cdots & 0 & 0 & 1 \\ 0 & \cdots & 0 & 1 & 0 \end{array}\right).
\]
Making the final connection and letting $\U = \IV{\{j_{2},\ldots,j_{n-2}\}}$, $w=\IV{j_{n-1}}$, $v=\IV{j_{n}}$ and $u=\IV{j_{1}}$, we see that we again have a  directed 3-cycle $(u,v,w,u)$ with bidirectional arcs between $u$ and $\U$, with $\U$ an independent set.  Thus, $R(d)\in \gT$, and the proof is complete.
\edprf



\section{Directed 3-cycle anchored degree sequences} \label{sec:anchored}

Up to this point, we have identified the ill-defined degree sequences for the parallel Havel-Hakimi algorithm as well as structurally their digraph realizations.  This section defines a new class of degree sequences, called $\DH$-anchored, and shows that the ill-defined degree sequences are precisely defined by this class, where $\DH$ is a directed 3-cycle $\CC$.

\begin{definition}
Given a degree sequence $d$ and digraph $\DH$, we say $d$ is {\bf potentially $\bm \DH$-digraphic} if and only if there exists $\DG\in R(d)$ with a subgraph $\DH^{\prime}\subseteq \DG$ such that $\DH^{\prime}\cong\DH$.  We say $d$ is {\bf forcibly $\bm \DH$-digraphic} if and only if this is satisfied for all $\DG\in R(d)$.
\end{definition}

We showed in Lemma~\ref{lem:gRgT} that every ill-defined degree sequence $d\in\gD$ is forcibly $\CC$-digraphic {\it locally} through the ill-defined indices.  It is a surprising fact that being forcibly $\CC$-digraphic locally through fixed indices is in fact sufficient for the degree sequence to be ill-defined with respect to $\HH^{\pm}$.  To prove this, we start with a definition.

\begin{definition}
Given a digraph $\DH$, we will call a degree sequence $d$ {\bf $\bm \DH$-anchored} if it is forcibly $\DH$-digraphic and there exists a nonempty set of indices $\mathcal{J}(\DH)$, called an {\bf $\bm \DH$-anchor set}, such that for every index $i \in \mathcal{J}(\DH)$ and every $\DG \in R(d)$, there is an induced subgraph $\DH^{\prime}\subseteq \DG$ with $\DH^{\prime}\cong \DH$ and $\IV{i} \in V(\DH^{\prime})$.
\label{def:anchor}
\end{definition}
We will also call a digraph $\DG$ $\DH$-anchored if $d_{\DG}$ is $\DH$-anchored.
\begin{theorem}
$d$ is $\CC$-anchored if and only if $d \in \mathcal{D}$.  More precisely, $i \in \mathcal{J}(\CC)$ if and only if $i$ is an ill-defined index for $\HH^{\pm}$.
\label{thm:anchor}
\end{theorem}
\bnprf
Let $d \in \mathcal{D}$ with $i$ an ill-defined index.  Since $R(d) \subseteq \gR$, by Lemma~\ref{lem:gRgT} there is a directed 3-cycle $\CC$ through $\IV{i}$ for every realization of $d$.  By definition, this implies $d$ is $\CC$-anchored with $i\in\mathcal{J}(\CC)$.

For the converse, suppose $d\notin\gD$ and let $i$ be a well-defined index.  Let $u = \IV{i}$ and $\Mp, \Mm$ maximal vertex sets for $u$. 
Since $i$ is well-defined, by Lemma~\ref{lem:gDgR}, $u$ is well-defined, so there is a $\DG^{\prime}\in R(d)$ such that $(\Mp,u,\Mm)\subseteq \DP(\DG^{\prime})$.  If there is a directed 3-cycle through $u$, it must be of the form $(h^{+},u,h^{-},h^{+})$ for some $h^{+}\in \Mp$ and $h^{-}\in \Mm$.  
Since $h^{+}\in \Mp$ and $h^{-}\notin \Mp$, we have $d^{+}_{h^{+}} \geq d^{+}_{h^{-}}$.  Thus, either $d^{+}_{h^{+}} = d^{+}_{h^{-}}$ with $d^{-}_{h^{+}} > d^{-}_{h^{-}}$, or $d^{+}_{h^{+}} > d^{+}_{h^{-}}$.  The first case can't happen since $h^{-} \in \Mm$ and $h^{+} \notin \Mm \Rightarrow d^{-}_{h^{-}} \geq d^{-}_{h^{+}}$, so we must have $d^{+}_{h^{+}} > d^{+}_{h^{-}}$.  Combining this with the fact that $(h^{-},h^{+})\in \DE$ means there is an $x \neq u$ such that $(h^{+},x)\in \DE$ and $(h^{-},x)\notin \DE$.  We now perform the 3-switch $\sigma_{3}\bigl((u,h^{-}),(h^{-},h^{+}),(h^{+},x)\bigr)$, thereby removing the directed 3-cycle.

This can be repeated for all directed 3-cycles through $u$, giving a realization of $d$ with no directed 3-cycles through $u$.  This shows $i\notin\mathcal{J}(\CC)$, and since $i$ was arbitrary, $d$ is not $\CC$-anchored.
\edprf

Given a degree sequence $d\in\gD$, from this point on we will call the set $\mathcal{J}(\CC)$ both the $\CC$-anchor set for $d$ and the set of ill-defined indices (relative to $\HH^{\pm}$).  The next few lemmas show that we can identify from the degree sequence which indices or vertices belong to which vertex class.  We start with a lemma that identifies the $\CC$-anchor set.

\begin{lem} Let $d$ be $\CC$-anchored with $\mathcal{J}(\CC)$ the $\CC$-anchor set.  Then $j \in \mathcal{J}(\CC)$ if and only if, in reference to the three cases (i)--(iii) in Definition~\ref{def:gS}, we have

(i) $j \in \{j_{1}, j_{2}, j_{3}\}$ and $d_{j} = (l+1,k+1)$

(ii) $j = j_{n}$ and $d_{j} = (l+1,k+1)$

(iii) $j = j_{1}$ and $d_{j} = (l+n-2,k+n-2)$.
\label{lem:illdef}
\end{lem}
\bnprf
Let $j$ be an ill-defined index for $d$ and $u=\IV{j}$.  Let $C=\{u,v,w\}\subseteq V(G)$ be a vertex set such that $\DG[C]\cong \CC$ for $\DG\in R(d)$.  Since $\gR=\gT$, we can decompose the vertex set $V(\DG)$ into the six vertex classes.  In the proof in Lemma~\ref{lem:gTgS}, we showed that $\U = \emptyset$ corresponds to case $(i)$ and $(ii)$ of Definition~\ref{def:gS} and in particular that $j = j_{n}$ with $d_{j} = (l+1,k+1)$.  If we are in case $(i)$, we have $d_{j_{1}} = d_{j_{2}} = d_{j_{3}}$, and so $j \in \{j_{1}, j_{2}, j_{3}\}$.  We omitted the case $\Z = \emptyset$ in Lemma~\ref{lem:gTgS}, but by a similar argument we arrive at case $(iii)$ of Definition~\ref{def:gS} with $j=j_{1}$ and $d_{j} = (l+n-2,k+n-2)$.

Conversely, by Lemma~\ref{lem:gSgR} we have that the indices in $(i)$--$(iii)$ are in fact ill-defined, and we have our result.
\edprf

\begin{definition}
Let $d$ be $\CC$-anchored with $\mathcal{J}(\CC)$ the $\CC$-anchor set.  For every $j\in\mathcal{J}(\CC)$, define the {\bf $\bm \CC$-scaffold set} $S_{j}$ by
\[
S_{j} = \bigcup_{\DG\in R(d)} \{\VI{C} : C \subseteq V(\DG) \ \mbox{with} \ \DG[C]\cong\CC \ \mbox{and} \ j\in\VI{C}\}.
\]
\end{definition}

\begin{lem}
Let $d$ be $\CC$-anchored with $\mathcal{J}(\CC)$ the $\CC$-anchor set.  Then for every $j \in \mathcal{J}(\CC)$, the $\CC$-scaffold set $S_{j}$ is given by $S_{j} = J$, where $J=\{j_{1},\ldots,j_{n}\}$ is the set of indices given by Definition~\ref{def:gS}.
\label{lem:scaffold}
\end{lem}
\bnprf
Let $j\in\mathcal{J}(\CC)$, $\DG = (V,\DE)\in R(d)$, $C=\{u,v,w\}\subseteq V$ such that $\DG[C]\cong\CC$ with $\VI{u}=j$, and $x\in V\setminus C$.  By Lemma~\ref{lem:illdef}, there is a set of indices $J = \{j_{1},\ldots,j_{n}\}$ satisfying Definition~\ref{def:gS} such that $j\in J$.  Since $\DG \in \gT$, from Lemma~\ref{lem:gTgS} we have that $J=\VI{C\cup\Z\cup\U}$.  Clearly $\VI{C}\subseteq S_{j}$, so suppose $\Z\neq\emptyset$.  Thus, we are in case $(ii)$ of Definition~\ref{def:gS}, and in Lemma~\ref{lem:gTgS} we showed $\VI{Z}=\{j_{1},j_{2}\}$, $\VI{\Z}=\{j_{3},\ldots,j_{n-1}\}$ and $\VI{u} = j \equiv j_{n}$.  By Table~\ref{tab:ineq}, $d_{\Z} = d_{v} = d_{w}$, and so for every $z\in\Z$ there is a $\DG^{\prime} \in R(d)$ such that $C^{\prime}=\{u,w,z\}$ with $\DG^{\prime}[C^{\prime}]\cong\CC$, which shows $\VI{\Z}\subseteq S_{j}$.  A similar argument holds for $\U$, and thus $J\subseteq S_{j}$.

For the converse, we need to be sure that for {\it all} directed 3-cycles $C$ containing $u$ that $\VI{C} \subseteq J$.  But the bottom diagrams in Figure~\ref{fig:3rules} show that all directed 3-cycles $C$ are such that $\VI{C} \subseteq J$ with $u\in C$, $\VI{C} \subseteq \VI{\Cp}$, or $\VI{C} \subseteq \VI{\Cm}$.  Thus we have $S_{j}\subseteq J$.
\edprf

The following corollary summarizes our results and shows how to identify the vertex classes in $\gR$ from the degree sequences in $\gD$.
\begin{cor}
Let $d$ be $\CC$-anchored with $\mathcal{J}(\CC)$ the $\CC$-anchor set.  Let $p$ and $q$ denote permutations such that $\bar{d}_{i} = d_{p_{i}}$ and $\ubar{d}_{i} = d_{q_{i}}$ and let $j \in \mathcal{J}(\CC)$.  Define the following sets
\begin{eqnarray*}
X_{1} & = & \{p_{1},\ldots,p_{k}\} \\
X_{2} & = & \{q_{1},\ldots,q_{l}\} \\
X_{3} & = & \{p_{k+n+1},\ldots,p_{N}\} \\
X_{4} & = & \{q_{l+n+1},\ldots,q_{N}\}.
\end{eqnarray*}
If $\DG_{d}\in R(d)$, $C=\{u,v,w\}\subseteq V(\DG)$ such that $\VI{u} = j$ and $\DG[C]\cong \CC$, then we have $\VI{u}$ given by Lemma~\ref{lem:illdef}, $\VI{Z\cup\Z\cup\U} = S_{j}-\{j\}$ by Lemma~\ref{lem:scaffold}, and by Eq.~(\ref{eq:classorder}), $\VI{\C} = X_{1}\cap X_{2}$, $\VI{\Cp} = X_{1}\setminus X_{2}$, $\VI{\Cm} = X_{2}\setminus X_{1}$, and $\VI{\I} = X_{3}\cap X_{4}$.
\label{cor:id_class}
\end{cor}



\section{The parallel algorithm revisited} \label{sec:algorithm}

The previous section completely classified both the degree sequences and their indices where the parallel Havel-Hakimi algorithm is ill-defined.  The ill-defined sequences are exactly the $\CC$-anchored degree sequences, with the $\CC$-anchor set corresponding to the ill-defined indices.  Through Definition \ref{def:gS}, we can identify the ill-defined sequences by their degree sequence characterizations, as well as the ill-defined indices in Lemma~\ref{lem:illdef}.  

At each step of the parallel algorithm, if the degree sequence has any ill-defined indices, we can proceed in two natural ways.  The first technique is to recognize that we know structurally what the realized digraph looks like {\it locally} around an ill-defined vertex $u$ from the definition of $\gT$, and since we can identify the vertex classes by Corollary \ref{cor:id_class}, we can connect them appropriately and {\it remove} the ill-defined index $i=\IV{u}$, which forces us to create a directed 3-cycle $\CC$ at $u$.  We then continue to remove ill-defined indices until the degree sequence itself again is well-defined or the residual degree sequence is the zero sequence.  This is shown in Algorithm~\ref{alg:alg1}.1.

\begin{figure}[h!]
\ttfamily\upshape

\begin{tabular}{l} \hline\hline
{\bf Algorithm~\ref{alg:alg1}.1:} Parallel Havel-Hakimi algorithm \\ \hline
{\bf Input:} Integer-pair sequence $d$ \\ 
{\bf Output:} Graph $\DG =(V,\DE)\in R(d)$ \\ \\
$N \coloneqq$ length of $d$ \\
$V \coloneqq \IV{\{1,\ldots,N\}}$ \\ 
$\DE \coloneqq \emptyset$ \\
{\bf if} {the slack sequence has a negative entry} {\bf then} \\
\qquad {\bf return} $\DG$ \\
{\bf end if} \\ \\
{\bf while} {$d$ is not the zero sequence} \\
\qquad {\bf if} {$\mathcal{J}(\CC) = \emptyset$} {\bf then}\\
\qquad\qquad $i \coloneqq$ index with non-zero degree \\
\qquad\qquad $\sK \coloneqq$ maximal index sets for $i$ \\
\qquad {\bf else} \\
\qquad\qquad $i \coloneqq$ index in $\mathcal{J}(\CC)$ \\
\qquad\qquad $\sK \coloneqq$ index sets constructed from Corollary~\ref{cor:id_class} \\
\qquad {\bf end if} \\
\qquad $[u,\sM] \coloneqq [\IV{i},\IV{\sK}]$ \\ 
\qquad $\DE \coloneqq \DE \cup (\Mp,u) \cup (u,\Mm)$ \\
\qquad $d \coloneqq \HH^{\pm}(d,i,\sK)$ \\
{\bf end while} \\ \\
{\bf return} $\DG$ \\
\hline\hline
\end{tabular}
\rmfamily

\label{alg:alg1}

\end{figure}

The second technique is to keep choosing well-defined indices.  However, we will at some point have to deal with the ill-defined indices since, as we will show below in Lemma~\ref{lem:persist}, the ill-defined indices persist once they appear in a residual degree sequence.  Once we have exhausted all well-defined indices, the resulting degree sequence will have only ill-defined indices, which we show to be unidigraphic in Lemma~\ref{lem:extseq} (a degree sequence is unidigraphic if it has only one isomorphism class in its set of digraph realizations).  Thus, when we arrive at one of these extreme degree sequences, the only ambiguity in finding a realization is in the orientation of the anchored directed 3-cycles.  This algorithm is displayed in Algorithm~\ref{alg:alg2}.2.

\begin{figure}[h!]
\ttfamily\upshape

\begin{tabular}{l} \hline\hline
{\bf Algorithm~\ref{alg:alg2}.2:} Alternate parallel Havel-Hakimi algorithm \\ \hline
{\bf Input:} Integer-pair sequence $d$ \\ 
{\bf Output:} Graph $\DG =(V,\DE)\in R(d)$ \\ \\
$N \coloneqq$ length of $d$ \\
$V \coloneqq \IV{\{1,\ldots,N\}}$ \\ 
$\DE \coloneqq \emptyset$ \\
{\bf if} {the slack sequence has a negative entry} {\bf then} \\
\qquad {\bf return} $\DG$ \\
{\bf end if} \\ \\
{\bf while} {$d$ is not the zero sequence} {\bf and} \\
\qquad \quad {$d$ has well-defined indices with non-zero degree} {\bf do}\\
\qquad $i \coloneqq$ well-defined index with non-zero degree \\
\qquad $\sK \coloneqq$ maximal index sets for $i$ \\
\qquad $[u,\sM] \coloneqq [\IV{i},\IV{\sK}]$ \\ 
\qquad $\DE \coloneqq \DE \cup (\Mp,u) \cup (u,\Mm)$ \\
\qquad $d \coloneqq \HH^{\pm}(d,i,\sK)$ \\
{\bf end while} \\ \\
{\bf if} {$\mathcal{J}(\CC) \neq \emptyset$} {\bf then} \\
\qquad $X \coloneqq$ all {\de}s given by extreme degree sequence (Lemma \ref{lem:extseq}) \\
\qquad $\DE \coloneqq \DE \cup X$ \\
{\bf end if} \\ \\
{\bf return} $\DG$ \\
\hline\hline
\end{tabular}
\rmfamily

\label{alg:alg2}
\end{figure}

\begin{lem} 
If $d$ is $\CC$-anchored, then for any well-defined index pair $[j,\sK]$, $\tilde{d} = \HH^{\pm}(d,j,\sK)$ is $\CC$-anchored with the $\CC$-anchor set preserved.
\label{lem:persist}
\end{lem}
\bnprf
Let $[j,\sK]$ be a well-defined index pair with $[x,\sM] = [\IV{j},\IV{\sK}]$ and let $i\in\mathcal{J}(\CC)$ with $u=\IV{i}$.  Since $[x,\sM]$ is a well-defined vertex pair, and by Lemmas~\ref{lem:illdef} and \ref{lem:scaffold}, there is a $\DG = (V,\DE) \in R(d)$ and $C=\{u,v,w\}$ such that $(\Mp,x,\Mm)\subseteq \DP(\DG)$, $\DG[C]\cong\CC$ and $x \in V\setminus C$.  One of the interesting properties of the $M$-partition on $V-C$ that follows immediately from the definition is that removing a vertex from $V-C$ preserves the $M$-partition structure.  Thus, we have $\tilde{d} = \HH^{\pm}(d,j,\sK)$ is $\CC$-anchored with the $\CC$-anchor set preserved.
\edprf

The following three digraphs are examples of realizations of the extreme degree sequences where all vertices are ill-defined, showing their recursive construction on the number of directed 3-cycles (note that $\CC$ is the simplest extreme case and is not shown).  The dashed circles with an arc between them means a unidirectional complete join.
\medskip
\begin{center}
\rotatebox{90}{\includegraphics[height=4in]{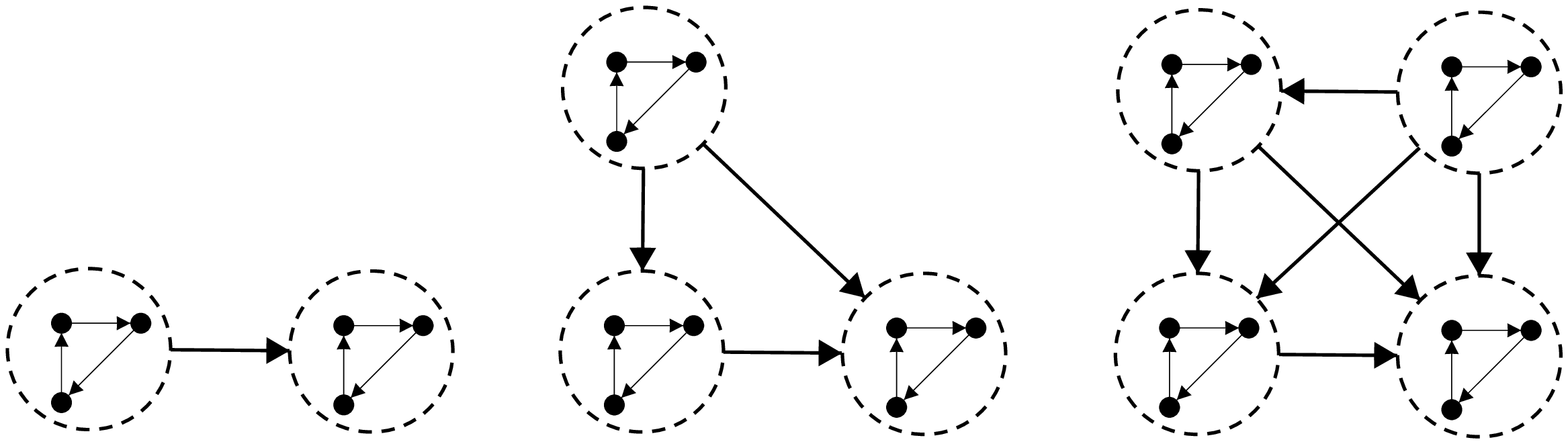}}
\end{center}
\medskip
\begin{lem}
The extreme degree sequences are unidigraphic.
\label{lem:extseq}
\end{lem}
\bnprf
Let $d$ be $\CC$-anchored such that $\mathcal{J}(\CC) = \VI{V}$, where $\DG = (V,\DE)\in R(d)$.  Since $\mathcal{J}(\CC) = \VI{V}$, we must be in case $(i)$ of Definition \ref{def:gS}, and thus by Lemma~\ref{lem:illdef} the number of vertices $N$ must be a multiple of 3.  Let $m$ be a positive integer such that $N = 3m$.  We have
\[
\begin{array}{c}
\bar{d}_{1} = \bar{d}_{2} = \bar{d}_{3} = (l_{1}, k_{1}) = (l_{1}, 1) \\
\bar{d}_{4} = \bar{d}_{5} = \bar{d}_{6} = (l_{2}, k_{2}) = (l_{2}, 4) \\
\vdots \\
\bar{d}_{3m-2} = \bar{d}_{3m-1} = \bar{d}_{3m} = (l_{m}, k_{m}) = (l_{m}, 3m-2).
\end{array}
\]
We also have
\[
\begin{array}{c}
\ubar{d}_{1} = \ubar{d}_{2} = \ubar{d}_{3} = (l_{m}, 3m-2) = (1, 3m-2) \\
\ubar{d}_{4} = \ubar{d}_{5} = \ubar{d}_{6} = (l_{m-1}, 3m-5) = (4, 3m-5) \\
\vdots \\
\ubar{d}_{3m-2} = \ubar{d}_{3m-1} = \ubar{d}_{3m} = (l_{1}, 1) = (3m-2, 1).
\end{array}
\]
We will use induction to prove that $d$ is unidigraphic.  For $m=1$, we have $d = \{(1,1),(1,1),(1,1)\}$, which is clearly unidigraphic.  For $m>1$, we know we have a directed 3-cycle $\CC$ between the vertices with degrees $\ubar{d}_{3m-2}$, $\ubar{d}_{3m-1}$ and $\ubar{d}_{3m}$.  After making those connections, we have
\[\ubar{d}_{3m-2} = \ubar{d}_{3m-1} = \ubar{d}_{3m} = (3(m-1),0).\]
Since there are exactly $3(m-1)$ indices left with $d^{-}$ non-zero, we connect their vertices and are left with
\[
\begin{array}{c}
\ubar{\tilde{d}}_{1} = \ubar{\tilde{d}}_{2} = \ubar{\tilde{d}}_{3} = (1, 3m-5) = (1, 3(m-1)-2) \\
\ubar{\tilde{d}}_{4} = \ubar{\tilde{d}}_{5} = \ubar{\tilde{d}}_{6} = (4, 3m-8) = (4, 3(m-1)-5) \\
\vdots \\
\ubar{\tilde{d}}_{3(m-1)-2} = \ubar{\tilde{d}}_{3(m-1)-1} = \ubar{\tilde{d}}_{3(m-1)} = (3m-5, 1) = (3(m-1)-2, 1).
\end{array}
\]
By induction, we have the result.
\edprf

As was seen in the proof of the theorem, the only ambiguity in how to connect vertices is in the orientation of the directed 3-cycles.



The original motivation for developing the parallel Havel-Hakimi algorithm is that for digraphic sequences where $\sum d^{+}$ is even and $d^{+}_{i}=d^{-}_{i}$ for all $i$, if at each step we choose maximal index sets $\sK$ such that $\Kp=\Km$, then the algorithm will realize a digraph $\DG$ with all bidirectional {\de}s.  In fact, in this special case the parallel Havel-Hakimi algorithm corresponds exactly with the undirected Havel-Hakimi algorithm if we identify bidirectional {\de}s with undirected edges.

\begin{definition}
Given an integer-pair sequence $d$, we say $d$ is {\bf Eulerian} if $d_{i}^{+} = d_{i}^{-}$ for all $i$.
\end{definition}

What digraph will the parallel Havel-Hakimi algorithm realize when a digraphic sequence $d$ is Eulerian, $\sum d^{+}$ is now {\it odd}, and we choose a maximal index pair at each step such that $\Kp = \Km$?  It turns out that the algorithm will realize a digraph with all bidirectional {\de}s {\it except} for one directed 3-cycle $\CC$, which we prove in Theorem~\ref{thm:eulerian} below.

\begin{theorem}
If a digraphic sequence $d$ is Eulerian with $\sum d^{+}$ odd, and at every step of the parallel Havel-Hakimi algorithm $\HH^{\pm}$ we choose maximal index sets $\sK$ such that $\Kp=\Km$, then the realization produced will have all bidirectional {\de}s except one 3-cycle.
\label{thm:eulerian}
\end{theorem}

\bnprf
Using Algorithm~\ref{alg:alg2}.2, if we always choose maximal index sets $\sK$ such that $\Kp=\Km$, then every step of $\HH^{\pm}$ creates bidirectional {\de}s and reduces the degree sum by an even number, with the residual degree sequence Eulerian as well.  Since the degree sum is reduced by an even number at every step and $\sum d^{+}$ is odd, at some point the residual degree sequence must be $\CC$-anchored.  Continue the algorithm until we arrive at an extreme degree sequence.  But by arguments in Lemma~\ref{lem:extseq} the only extreme degree sequence that is Eulerian is the degree sequence for a single directed 3-cycle, and thus we have the theorem.
\edprf

\section*{Acknowledgments}
The author thanks Gexin Yu and Laura Sheppardson for helpful discussions.  This work was funded by the postdoctoral and undergraduate biological sciences education program grant awarded to the College of William and Mary by the Howard Hughes Medical Institute.

\appendix
\section*{Appendix: Seeds for Induction}
\label{appendix}

Table~\ref{tab:R1} includes the 16 cases for the proof in Lemma~\ref{lem:seed} showing that $\gQ_1 \subseteq \gT_1$, with each case describing how a vertex $x$ relates to a directed 3-cycle $\CC$ (see Lemma~\ref{lem:seed} for a description of the cases).  The second column lists the degree sequence in matrix form, where the columns correspond to $(d_{u},d_{v},d_{w},d_{x})$.  A checkmark ($\checkmark$) in the third column implies the case is in $\gQ_1$ with the corresponding ill-defined maximal vertex sets given in the last column.  A dot ($\cdot$) in the third column implies the case is not in $\gQ_1$, with the last column stating which part of Lemma~\ref{lem:gQ} is violated for every pair of maximal sets.

Table~\ref{tab:R2} includes the 84 cases for the proof in Lemma~\ref{lem:seed} showing that $\gQ_2 \subseteq \gT_2$, with each case denoted in the first column by $[X,Y]$ such that $X,Y\in\{\I,\U,\Z,\Cp,\Cm,\C\}$.  For each case in the first column, the second column lists the 4 possible relations between a vertex $x \in X$ and $y\in Y$. The third column lists the degree sequence in matrix form, where the columns correspond to $(d_{u},d_{v},d_{w},d_{x},d_{y})$.  A checkmark ($\checkmark$) in the fourth column implies the case is in $\gQ_2$ with the corresponding ill-defined maximal vertex sets given in the last column.  A dot ($\cdot$) in the fourth column implies the case is not in $\gQ_2$, with the last column stating which part of Lemma~\ref{lem:gQ} is violated for every pair of maximal sets.  See Figure~\ref{fig:gTable} for a graphic summary of this table.

\begin{table}[t]
\[
\begin{array}{|c|c|c|l|}
\hline
(U^{0},Z^{0}) & {\tiny \left(\begin{array}{c} 1 \ 1 \ 1 \ 0 \\ 1 \ 1 \ 1 \ 0 \end{array}\right)} & \checkmark & \Mp = \Mm = \{w\} \\ \hline
(U^{0},Z^{+}) & {\tiny \left(\begin{array}{c} 1 \ 1 \ 1 \ 2 \\ 1 \ 2 \ 2 \ 0 \end{array}\right)} & \cdot & (ii) \ w\notin \Mp\cap\Mm \\ \hline
(U^{0},Z^{-}) & {\tiny \left(\begin{array}{c} 1 \ 2 \ 2 \ 0 \\ 1 \ 1 \ 1 \ 2 \end{array}\right)} & \cdot & \ooalign{\hfil\phantom{$(ii) \ w\notin \Mp\cap\Mm$}\hfil\cr \hfil$''$\hfil\cr} \\ \hline
(U^{0},Z^{\pm}) & {\tiny \left(\begin{array}{c} 1 \ 2 \ 2 \ 2 \\ 1 \ 2 \ 2 \ 2 \end{array}\right)} & \checkmark & \Mp = \Mm = \{w\} \\ \hline\hline

(U^{+},Z^{0}) & {\tiny \left(\begin{array}{c} 1 \ 1 \ 1 \ 1 \\ 2 \ 1 \ 1 \ 0 \end{array}\right)} & \cdot & (ii) \ v \in \Mp \\ \hline
(U^{+},Z^{+}) & {\tiny \left(\begin{array}{c} 1 \ 1 \ 1 \ 3 \\ 2 \ 2 \ 2 \ 0 \end{array}\right)} & \checkmark & \Mp=\{w,x\}, \ \Mm=\{w\} \\ \hline
(U^{+},Z^{-}) & {\tiny \left(\begin{array}{c} 1 \ 2 \ 2 \ 1 \\ 2 \ 1 \ 1 \ 2 \end{array}\right)} & \cdot & (ii) \ w\notin \Mp\cap\Mm \\ \hline
(U^{+},Z^{\pm}) & {\tiny \left(\begin{array}{c} 1 \ 2 \ 2 \ 3 \\ 2 \ 2 \ 2 \ 2 \end{array}\right)} & \cdot & (ii) \ w\notin \Mp\cap\Mm \\ \hline\hline

(U^{-},Z^{0}) & {\tiny \left(\begin{array}{c} 2 \ 1 \ 1 \ 0 \\ 1 \ 1 \ 1 \ 1 \end{array}\right)} & \cdot & (ii) \ v \in \Mm \\ \hline
(U^{-},Z^{+}) & {\tiny \left(\begin{array}{c} 2 \ 1 \ 1 \ 2 \\ 1 \ 2 \ 2 \ 1 \end{array}\right)} & \cdot & (ii) \ w\notin \Mp\cap\Mm \\ \hline
(U^{-},Z^{-}) & {\tiny \left(\begin{array}{c} 2 \ 2 \ 2 \ 0 \\ 1 \ 1 \ 1 \ 3 \end{array}\right)} & \checkmark & \Mp=\{w\}, \ \Mm=\{w,x\} \\ \hline
(U^{-},Z^{\pm}) & {\tiny \left(\begin{array}{c} 2 \ 2 \ 2 \ 2 \\ 1 \ 2 \ 2 \ 3 \end{array}\right)} & \cdot & (ii) \ w\notin \Mp\cap\Mm \\ \hline\hline

(U^{\pm},Z^{0}) & {\tiny \left(\begin{array}{c} 2 \ 1 \ 1 \ 1 \\ 2 \ 1 \ 1 \ 1 \end{array}\right)} & \checkmark & \Mp = \Mm = \{w\} \\ \hline
(U^{\pm},Z^{+}) & {\tiny \left(\begin{array}{c} 2 \ 1 \ 1 \ 3 \\ 2 \ 2 \ 2 \ 1 \end{array}\right)} & \cdot & (ii) \ v \in \Mm \\ \hline
(U^{\pm},Z^{-}) & {\tiny \left(\begin{array}{c} 2 \ 2 \ 2 \ 1 \\ 2 \ 1 \ 1 \ 3 \end{array}\right)} & \cdot & (ii) \ v \in \Mp \\ \hline
(U^{\pm},Z^{\pm}) & {\tiny \left(\begin{array}{c} 2 \ 2 \ 2 \ 3 \\ 2 \ 2 \ 2 \ 3 \end{array}\right)} & \checkmark & \Mp = \Mm = \{w\} \\ \hline
\end{array}
\]
\caption{The 16 cases for the proof in Lemma~\ref{lem:seed} showing that $\gQ_1 \subseteq \gT_1$, with each case describing how a vertex $x$ relates to a directed 3-cycle $\CC$ (see Lemma~\ref{lem:seed} for a description of the cases).  A description of the table is included in the text for Appendix~\ref{appendix}.}
\label{tab:R1}
\end{table}

\begin{table}[t]
\[
\begin{array}{|c|c|c|c|l|}
\hline
\multirow{4}{*}{$[\I,\I]$} & x \cdots y & {\tiny \left(\begin{array}{c} 1 \ 1 \ 1 \ 0 \ 0 \\ 1 \ 1 \ 1 \ 0 \ 0 \end{array}\right)} & \checkmark & \Mp = \Mm = \{w\} \\ \cline{2-5}
 & x \rightarrow y & {\tiny \left(\begin{array}{c} 1 \ 1 \ 1 \ 1 \ 0 \\ 1 \ 1 \ 1 \ 0 \ 1 \end{array}\right)} & \cdot & (v) \ (x,v)\notin \DE \ \mbox{and} \ (v,y)\notin \DE\\ \cline{2-5}
 & x \leftarrow y & {\tiny \left(\begin{array}{c} 1 \ 1 \ 1 \ 0 \ 1 \\ 1 \ 1 \ 1 \ 1 \ 0 \end{array}\right)} & \cdot & (v) \ (x,v)\notin \DE \ \mbox{and} \ (v,y)\notin \DE \\ \cline{2-5}
 & x \leftrightarrow y & {\tiny \left(\begin{array}{c} 1 \ 1 \ 1 \ 1 \ 1 \\ 1 \ 1 \ 1 \ 1 \ 1 \end{array}\right)} & \cdot & (v) \ (x,v)\notin \DE \ \mbox{and} \ (v,y)\notin \DE \\ \hline\hline

\multirow{4}{*}{$[\U,\I]$} & x \cdots y & {\tiny \left(\begin{array}{c} 2 \ 1 \ 1 \ 1 \ 0 \\ 2 \ 1 \ 1 \ 1 \ 0 \end{array}\right)} & \checkmark & \Mp = \Mm = \{w,x\} \\ \cline{2-5}
 & x \rightarrow y & {\tiny \left(\begin{array}{c} 2 \ 1 \ 1 \ 2 \ 0 \\ 2 \ 1 \ 1 \ 1 \ 1 \end{array}\right)} & \cdot & (v) \ (x,v)\notin \DE \ \mbox{and} \ (v,y)\notin \DE \\ \cline{2-5}
 & x \leftarrow y & {\tiny \left(\begin{array}{c} 2 \ 1 \ 1 \ 1 \ 1 \\ 2 \ 1 \ 1 \ 2 \ 0 \end{array}\right)} & \cdot & (v) \ (x,v)\notin \DE \ \mbox{and} \ (v,y)\notin \DE \\ \cline{2-5}
 & x \leftrightarrow y & {\tiny \left(\begin{array}{c} 2 \ 1 \ 1 \ 2 \ 1 \\ 2 \ 1 \ 1 \ 2 \ 1 \end{array}\right)} & \cdot & (v) \ (x,v)\notin \DE \ \mbox{and} \ (v,y)\notin \DE \\ \hline\hline

\multirow{4}{*}{$[\U,\U]$} & x \cdots y & {\tiny \left(\begin{array}{c} 3 \ 1 \ 1 \ 1 \ 1 \\ 3 \ 1 \ 1 \ 1 \ 1 \end{array}\right)} & \checkmark & \Mp = \Mm = \{w,x,y\} \\ \cline{2-5}
 & x \rightarrow y & {\tiny \left(\begin{array}{c} 3 \ 1 \ 1 \ 2 \ 1 \\ 3 \ 1 \ 1 \ 1 \ 2 \end{array}\right)} & \cdot & (v) \ (x,v)\notin \DE \ \mbox{and} \ (v,y)\notin \DE \\ \cline{2-5}
 & x \leftarrow y & {\tiny \left(\begin{array}{c} 3 \ 1 \ 1 \ 1 \ 2 \\ 3 \ 1 \ 1 \ 2 \ 1 \end{array}\right)} & \cdot & (v) \ (x,v)\notin \DE \ \mbox{and} \ (v,y)\notin \DE \\ \cline{2-5}
 & x \leftrightarrow y & {\tiny \left(\begin{array}{c} 3 \ 1 \ 1 \ 2 \ 2 \\ 3 \ 1 \ 1 \ 2 \ 2 \end{array}\right)} & \cdot & (v) \ (x,v)\notin \DE \ \mbox{and} \ (v,y)\notin \DE \\ \hline\hline

\multirow{4}{*}{$[\Z,\I]$} & x \cdots y & {\tiny \left(\begin{array}{c} 1 \ 2 \ 2 \ 2 \ 0 \\ 1 \ 2 \ 2 \ 2 \ 0 \end{array}\right)} & \checkmark & \Mp = \Mm = \{w\} \\ \cline{2-5}
 & x \rightarrow y & {\tiny \left(\begin{array}{c} 1 \ 2 \ 2 \ 3 \ 0 \\ 1 \ 2 \ 2 \ 2 \ 1 \end{array}\right)} & \cdot & (ii) \ w\notin \Mp\cap\Mm \\ \cline{2-5}
 & x \leftarrow y & {\tiny \left(\begin{array}{c} 1 \ 2 \ 2 \ 2 \ 1 \\ 1 \ 2 \ 2 \ 3 \ 0 \end{array}\right)} & \cdot & (ii) \ w\notin \Mp\cap\Mm \\ \cline{2-5}
 & x \leftrightarrow y & {\tiny \left(\begin{array}{c} 1 \ 2 \ 2 \ 3 \ 1 \\ 1 \ 2 \ 2 \ 3 \ 1 \end{array}\right)} & \cdot & (ii) \ w\notin \Mp\cap\Mm \\ \hline\hline

\multirow{4}{*}{$[\Z,\U]$} & x \cdots y & {\tiny \left(\begin{array}{c} 2 \ 2 \ 2 \ 2 \ 1 \\ 2 \ 2 \ 2 \ 2 \ 1 \end{array}\right)} & \cdot & (vi) \ (u,y)\in\DE \ \mbox{and} \ y\notin \Mm\\ \cline{2-5}
 & x \rightarrow y & {\tiny \left(\begin{array}{c} 2 \ 2 \ 2 \ 3 \ 1 \\ 2 \ 2 \ 2 \ 2 \ 2 \end{array}\right)} & \cdot & (vi) \ (u,y)\in\DE \ \mbox{and} \ y\notin \Mm \\ \cline{2-5}
 & x \leftarrow y & {\tiny \left(\begin{array}{c} 2 \ 2 \ 2 \ 2 \ 2 \\ 2 \ 2 \ 2 \ 3 \ 1 \end{array}\right)} & \cdot & (vi) \ (u,y)\in\DE \ \mbox{and} \ y\notin \Mm \\ \cline{2-5}
 & x \leftrightarrow y & {\tiny \left(\begin{array}{c} 2 \ 2 \ 2 \ 3 \ 2 \\ 2 \ 2 \ 2 \ 3 \ 2 \end{array}\right)} & \cdot & (vi) \ x\in\Mm \ \mbox{and} \ (u,x)\notin \DE \\ \hline\hline

\multirow{4}{*}{$[\Z,\Z]$} & x \cdots y & {\tiny \left(\begin{array}{c} 1 \ 3 \ 3 \ 2 \ 2 \\ 1 \ 3 \ 3 \ 2 \ 2 \end{array}\right)} & \cdot & (vii) \ (y,w)\in\DE \ \mbox{and} \ (y,x)\notin\DE\\ \cline{2-5}
 & x \rightarrow y & {\tiny \left(\begin{array}{c} 1 \ 3 \ 3 \ 3 \ 2 \\ 1 \ 3 \ 3 \ 2 \ 3 \end{array}\right)} & \cdot & (vii) \ (y,w)\in\DE \ \mbox{and} \ (y,x)\notin\DE \\ \cline{2-5}
 & x \leftarrow y & {\tiny \left(\begin{array}{c} 1 \ 3 \ 3 \ 2 \ 3 \\ 1 \ 3 \ 3 \ 3 \ 2 \end{array}\right)} & \cdot & (vii) \ (w,y)\in\DE \ \mbox{and} \ (x,y)\notin\DE \\ \cline{2-5}
 & x \leftrightarrow y & {\tiny \left(\begin{array}{c} 1 \ 3 \ 3 \ 3 \ 3 \\ 1 \ 3 \ 3 \ 3 \ 3 \end{array}\right)} & \checkmark & \Mp=\Mm=\{w\} \\ \hline\hline

\multirow{4}{*}{$[\Cm,\I]$} & x \cdots y & {\tiny \left(\begin{array}{c} 2 \ 2 \ 2 \ 0 \ 0 \\ 1 \ 1 \ 1 \ 3 \ 0 \end{array}\right)} & \checkmark & \Mp=\{w\},\ \Mm=\{w,x\}\\ \cline{2-5}
 & x \rightarrow y & {\tiny \left(\begin{array}{c} 2 \ 2 \ 2 \ 1 \ 0 \\ 1 \ 1 \ 1 \ 3 \ 1 \end{array}\right)} & \cdot & (v) \ (x,v)\notin \DE \ \mbox{and} \ (v,y)\notin \DE \\ \cline{2-5}
 & x \leftarrow y & {\tiny \left(\begin{array}{c} 2 \ 2 \ 2 \ 0 \ 1 \\ 1 \ 1 \ 1 \ 4 \ 0 \end{array}\right)} & \checkmark & \Mp=\{w\},\ \Mm=\{w,x\} \\ \cline{2-5}
 & x \leftrightarrow y & {\tiny \left(\begin{array}{c} 2 \ 2 \ 2 \ 1 \ 1 \\ 1 \ 1 \ 1 \ 4 \ 1 \end{array}\right)} & \cdot & (v) \ (x,v)\notin \DE \ \mbox{and} \ (v,y)\notin \DE \\ \hline\hline

\multirow{4}{*}{$[\Cm,\U]$} & x \cdots y & {\tiny \left(\begin{array}{c} 3 \ 2 \ 2 \ 0 \ 1 \\ 2 \ 1 \ 1 \ 3 \ 1 \end{array}\right)} & \cdot & (ii) \ v\in\Mp \\ \cline{2-5}
 & x \rightarrow y & {\tiny \left(\begin{array}{c} 3 \ 2 \ 2 \ 1 \ 1 \\ 2 \ 1 \ 1 \ 3 \ 2 \end{array}\right)} & \cdot & (ii) \ v\in\Mp \\ \cline{2-5}
 & x \leftarrow y & {\tiny \left(\begin{array}{c} 3 \ 2 \ 2 \ 0 \ 2 \\ 2 \ 1 \ 1 \ 4 \ 1 \end{array}\right)} & \checkmark & \Mp=\{w,y\},\ \Mm=\{w,x,y\} \\ \cline{2-5}
 & x \leftrightarrow y & {\tiny \left(\begin{array}{c} 3 \ 2 \ 2 \ 1 \ 2 \\ 2 \ 1 \ 1 \ 4 \ 2 \end{array}\right)} & \cdot & (vii) \ (x,y)\in\DE \ \mbox{and} \ (x,w)\notin\DE \\ \hline\hline

\multirow{4}{*}{$[\Cm,\Z]$} & x \cdots y & {\tiny \left(\begin{array}{c} 2 \ 3 \ 3 \ 0 \ 2 \\ 1 \ 2 \ 2 \ 3 \ 2 \end{array}\right)} & \cdot & (vii) \ (w,x)\in\DE \ \mbox{and} \ (y,x)\notin\DE \\ \cline{2-5}
 & x \rightarrow y & {\tiny \left(\begin{array}{c} 2 \ 3 \ 3 \ 1 \ 2 \\ 1 \ 2 \ 2 \ 3 \ 3 \end{array}\right)} & \cdot & (ii) \ w\notin\Mp\cap\Mm \\ \cline{2-5}
 & x \leftarrow y & {\tiny \left(\begin{array}{c} 2 \ 3 \ 3 \ 0 \ 3 \\ 1 \ 2 \ 2 \ 4 \ 2 \end{array}\right)} & \checkmark & \Mp=\{w\},\ \Mm=\{w,x\} \\ \cline{2-5}
 & x \leftrightarrow y & {\tiny \left(\begin{array}{c} 2 \ 3 \ 3 \ 1 \ 3 \\ 1 \ 2 \ 2 \ 4 \ 3 \end{array}\right)} & \cdot & (ii) \ w\notin\Mp\cap\Mm \\ \hline\hline

\multirow{4}{*}{$[\Cm,\Cm]$} & x \cdots y & {\tiny \left(\begin{array}{c} 3 \ 3 \ 3 \ 0 \ 0 \\ 1 \ 1 \ 1 \ 3 \ 3 \end{array}\right)} & \checkmark & \Mp=\{w\},\ \Mm=\{w,x,y\} \\ \cline{2-5}
 & x \rightarrow y & {\tiny \left(\begin{array}{c} 3 \ 3 \ 3 \ 1 \ 0 \\ 1 \ 1 \ 1 \ 3 \ 4 \end{array}\right)} & \checkmark & \Mp=\{w\},\ \Mm=\{w,x,y\} \\ \cline{2-5}
 & x \leftarrow y & {\tiny \left(\begin{array}{c} 3 \ 3 \ 3 \ 0 \ 1 \\ 1 \ 1 \ 1 \ 4 \ 3 \end{array}\right)} & \checkmark & \Mp=\{w\},\ \Mm=\{w,x,y\} \\ \cline{2-5}
 & x \leftrightarrow y & {\tiny \left(\begin{array}{c} 3 \ 3 \ 3 \ 1 \ 1 \\ 1 \ 1 \ 1 \ 4 \ 4 \end{array}\right)} & \checkmark & \Mp=\{w\},\ \Mm=\{w,x,y\} \\ \hline
\end{array}
\]
\caption{The 84 cases for the proof in Lemma~\ref{lem:seed} showing that $\gQ_2 \subseteq \gT_2$, with each case denoted in the first column by $[X,Y]$ such that $X,Y\in\{\I,\U,\Z,\Cp,\Cm,\C\}$.  A description of the table is included in the text for Appendix~\ref{appendix}.}
\label{tab:R2}
\end{table}

\begin{table}[t]
\[
\begin{array}{|c|c|c|c|l|}
\hline
\multirow{4}{*}{$[\Cp,\I]$} & x \cdots y & {\tiny \left(\begin{array}{c} 1 \ 1 \ 1 \ 3 \ 0 \\ 2 \ 2 \ 2 \ 0 \ 0 \end{array}\right)} & \checkmark & \Mp=\{w,x\},\ \Mm=\{w\} \\ \cline{2-5}
 & x \rightarrow y & {\tiny \left(\begin{array}{c} 1 \ 1 \ 1 \ 4 \ 0 \\ 2 \ 2 \ 2 \ 0 \ 1 \end{array}\right)} & \checkmark & \Mp=\{w,x\},\ \Mm=\{w\} \\ \cline{2-5}
 & x \leftarrow y & {\tiny \left(\begin{array}{c} 1 \ 1 \ 1 \ 3 \ 1 \\ 2 \ 2 \ 2 \ 1 \ 0 \end{array}\right)} & \cdot & (v) \ (y,v)\notin\DE \ \mbox{and} \ (v,x)\notin\DE \\ \cline{2-5}
 & x \leftrightarrow y & {\tiny \left(\begin{array}{c} 1 \ 1 \ 1 \ 4 \ 1 \\ 2 \ 2 \ 2 \ 1 \ 1 \end{array}\right)} & \cdot & (v) \ (y,v)\notin\DE \ \mbox{and} \ (v,x)\notin\DE \\ \hline\hline

\multirow{4}{*}{$[\Cp,\U]$} & x \cdots y & {\tiny \left(\begin{array}{c} 2 \ 1 \ 1 \ 3 \ 1 \\ 3 \ 2 \ 2 \ 0 \ 1 \end{array}\right)} & \cdot & (ii) \ v\in\Mm \\ \cline{2-5}
 & x \rightarrow y & {\tiny \left(\begin{array}{c} 2 \ 1 \ 1 \ 4 \ 1 \\ 3 \ 2 \ 2 \ 0 \ 2 \end{array}\right)} & \checkmark & \Mp=\{w,x,y\},\ \Mm=\{w,y\} \\ \cline{2-5}
 & x \leftarrow y & {\tiny \left(\begin{array}{c} 2 \ 1 \ 1 \ 3 \ 2 \\ 3 \ 2 \ 2 \ 1 \ 1 \end{array}\right)} & \cdot & (ii) \ v\in\Mm \\ \cline{2-5}
 & x \leftrightarrow y & {\tiny \left(\begin{array}{c} 2 \ 1 \ 1 \ 4 \ 2 \\ 3 \ 2 \ 2 \ 1 \ 2 \end{array}\right)} & \cdot & (v) \ (y,v)\notin\DE \ \mbox{and} \ (v,x)\notin\DE \\ \hline\hline

\multirow{4}{*}{$[\Cp,\Z]$} & x \cdots y & {\tiny \left(\begin{array}{c} 1 \ 2 \ 2 \ 3 \ 2 \\ 2 \ 3 \ 3 \ 0 \ 2 \end{array}\right)} & \cdot & (vii) \ (x,w)\in\DE \ \mbox{and} \ (x,y)\notin\DE \\ \cline{2-5}
 & x \rightarrow y & {\tiny \left(\begin{array}{c} 1 \ 2 \ 2 \ 4 \ 2 \\ 2 \ 3 \ 3 \ 0 \ 3 \end{array}\right)} & \checkmark & \Mp=\{w,x\},\ \Mm=\{w\} \\ \cline{2-5}
 & x \leftarrow y & {\tiny \left(\begin{array}{c} 1 \ 2 \ 2 \ 3 \ 3 \\ 2 \ 3 \ 3 \ 1 \ 2 \end{array}\right)} & \cdot & (ii) \ w\notin\Mp\cap\Mm \\ \cline{2-5}
 & x \leftrightarrow y & {\tiny \left(\begin{array}{c} 1 \ 2 \ 2 \ 4 \ 3 \\ 2 \ 3 \ 3 \ 1 \ 3 \end{array}\right)} & \cdot & (ii) \ w\notin\Mp\cap\Mm \\ \hline\hline

\multirow{4}{*}{$[\Cp,\Cm]$} & x \cdots y & {\tiny \left(\begin{array}{c} 2 \ 2 \ 2 \ 3 \ 0 \\ 2 \ 2 \ 2 \ 0 \ 3 \end{array}\right)} & \cdot & (vii) \ (x,w)\in\DE \ \mbox{and} \ (x,y)\notin\DE \\ \cline{2-5}
 & x \rightarrow y & {\tiny \left(\begin{array}{c} 2 \ 2 \ 2 \ 4 \ 0 \\ 2 \ 2 \ 2 \ 0 \ 4 \end{array}\right)} & \checkmark & \Mp=\{w,x\},\ \Mm=\{w,y\} \\ \cline{2-5}
 & x \leftarrow y & {\tiny \left(\begin{array}{c} 2 \ 2 \ 2 \ 3 \ 1 \\ 2 \ 2 \ 2 \ 1 \ 3 \end{array}\right)} & \cdot & (vii) \ (x,w)\in\DE \ \mbox{and} \ (x,y)\notin\DE \\ \cline{2-5}
 & x \leftrightarrow y & {\tiny \left(\begin{array}{c} 2 \ 2 \ 2 \ 4 \ 1 \\ 2 \ 2 \ 2 \ 1 \ 4 \end{array}\right)} & \cdot & (v) \ (y,v)\notin\DE \ \mbox{and} \ (v,x)\notin\DE \\ \hline\hline

\multirow{4}{*}{$[\Cp,\Cp]$} & x \cdots y & {\tiny \left(\begin{array}{c} 1 \ 1 \ 1 \ 3 \ 3 \\ 3 \ 3 \ 3 \ 0 \ 0 \end{array}\right)} & \checkmark & \Mp=\{w,x,y\},\ \Mm=\{w\} \\ \cline{2-5}
 & x \rightarrow y & {\tiny \left(\begin{array}{c} 1 \ 1 \ 1 \ 4 \ 3 \\ 3 \ 3 \ 3 \ 0 \ 1 \end{array}\right)} & \checkmark & \Mp=\{w,x,y\},\ \Mm=\{w\} \\ \cline{2-5}
 & x \leftarrow y & {\tiny \left(\begin{array}{c} 1 \ 1 \ 1 \ 3 \ 4 \\ 3 \ 3 \ 3 \ 1 \ 0 \end{array}\right)} & \checkmark & \Mp=\{w,x,y\},\ \Mm=\{w\} \\ \cline{2-5}
 & x \leftrightarrow y & {\tiny \left(\begin{array}{c} 1 \ 1 \ 1 \ 4 \ 4 \\ 3 \ 3 \ 3 \ 1 \ 1 \end{array}\right)} & \checkmark & \Mp=\{w,x,y\},\ \Mm=\{w\} \\ \hline\hline

\multirow{4}{*}{$[\C,\I]$} & x \cdots y & {\tiny \left(\begin{array}{c} 2 \ 2 \ 2 \ 3 \ 0 \\ 2 \ 2 \ 2 \ 3 \ 0 \end{array}\right)} & \checkmark & \Mp=\Mm=\{w,x\} \\ \cline{2-5}
 & x \rightarrow y & {\tiny \left(\begin{array}{c} 2 \ 2 \ 2 \ 4 \ 0 \\ 2 \ 2 \ 2 \ 3 \ 1 \end{array}\right)} & \checkmark & \Mp=\Mm=\{w,x\} \\ \cline{2-5}
 & x \leftarrow y & {\tiny \left(\begin{array}{c} 2 \ 2 \ 2 \ 3 \ 1 \\ 2 \ 2 \ 2 \ 4 \ 0 \end{array}\right)} & \checkmark & \Mp=\Mm=\{w,x\} \\ \cline{2-5}
 & x \leftrightarrow y & {\tiny \left(\begin{array}{c} 2 \ 2 \ 2 \ 4 \ 1 \\ 2 \ 2 \ 2 \ 4 \ 1 \end{array}\right)} & \checkmark & \Mp=\Mm=\{w,x\} \\ \hline\hline

\multirow{4}{*}{$[\C,\U]$} & x \cdots y & {\tiny \left(\begin{array}{c} 3 \ 2 \ 2 \ 3 \ 1 \\ 3 \ 2 \ 2 \ 3 \ 1 \end{array}\right)} & \cdot & (ii) \ v\in\Mp \\ \cline{2-5}
 & x \rightarrow y & {\tiny \left(\begin{array}{c} 3 \ 2 \ 2 \ 4 \ 1 \\ 3 \ 2 \ 2 \ 3 \ 2 \end{array}\right)} & \cdot & (ii) \ v\in\Mp \\ \cline{2-5}
 & x \leftarrow y & {\tiny \left(\begin{array}{c} 3 \ 2 \ 2 \ 3 \ 2 \\ 3 \ 2 \ 2 \ 4 \ 1 \end{array}\right)} & \cdot & (ii) \ v\in\Mp \\ \cline{2-5}
 & x \leftrightarrow y & {\tiny \left(\begin{array}{c} 3 \ 2 \ 2 \ 4 \ 2 \\ 3 \ 2 \ 2 \ 4 \ 2 \end{array}\right)} & \checkmark & \Mp=\Mm=\{w,x,y\} \\ \hline\hline

\multirow{4}{*}{$[\C,\Z]$} & x \cdots y & {\tiny \left(\begin{array}{c} 2 \ 3 \ 3 \ 3 \ 2 \\ 2 \ 3 \ 3 \ 3 \ 2 \end{array}\right)} & \cdot & (vii) \ (w,x)\in\DE \ \mbox{and} \ (y,x)\notin\DE \\ \cline{2-5}
 & x \rightarrow y & {\tiny \left(\begin{array}{c} 2 \ 3 \ 3 \ 4 \ 2 \\ 2 \ 3 \ 3 \ 3 \ 3 \end{array}\right)} & \cdot & (vii) \ (w,x)\in\DE \ \mbox{and} \ (y,x)\notin\DE \\ \cline{2-5}
 & x \leftarrow y & {\tiny \left(\begin{array}{c} 2 \ 3 \ 3 \ 3 \ 3 \\ 2 \ 3 \ 3 \ 4 \ 2 \end{array}\right)} & \cdot & (vii) \ (x,w)\in\DE \ \mbox{and} \ (x,y)\notin\DE \\ \cline{2-5}
 & x \leftrightarrow y & {\tiny \left(\begin{array}{c} 2 \ 3 \ 3 \ 4 \ 3 \\ 2 \ 3 \ 3 \ 4 \ 3 \end{array}\right)} & \checkmark & \Mp=\Mm=\{w,x\} \\ \hline\hline

\multirow{4}{*}{$[\C,\Cm]$} & x \cdots y & {\tiny \left(\begin{array}{c} 3 \ 3 \ 3 \ 3 \ 0 \\ 2 \ 2 \ 2 \ 3 \ 3 \end{array}\right)} & \cdot & (vii) \ (w,y)\in\DE \ \mbox{and} \ (x,y)\notin\DE \\ \cline{2-5}
 & x \rightarrow y & {\tiny \left(\begin{array}{c} 3 \ 3 \ 3 \ 4 \ 0 \\ 2 \ 2 \ 2 \ 3 \ 4 \end{array}\right)} & \checkmark & \Mp=\{w,x\},\ \Mm=\{w,x,y\} \\ \cline{2-5}
 & x \leftarrow y & {\tiny \left(\begin{array}{c} 3 \ 3 \ 3 \ 3 \ 1 \\ 2 \ 2 \ 2 \ 4 \ 3 \end{array}\right)} & \cdot & (vii) \ (w,y)\in\DE \ \mbox{and} \ (x,y)\notin\DE \\ \cline{2-5}
 & x \leftrightarrow y & {\tiny \left(\begin{array}{c} 3 \ 3 \ 3 \ 4 \ 1 \\ 2 \ 2 \ 2 \ 4 \ 4 \end{array}\right)} & \checkmark & \Mp=\{w,x\},\ \Mm=\{w,x,y\} \\ \hline\hline

\multirow{4}{*}{$[\C,\Cp]$} & x \cdots y & {\tiny \left(\begin{array}{c} 2 \ 2 \ 2 \ 3 \ 3 \\ 3 \ 3 \ 3 \ 3 \ 0 \end{array}\right)} & \cdot & (vii) \ (y,w)\in\DE \ \mbox{and} \ (y,x)\notin\DE \\ \cline{2-5}
 & x \rightarrow y & {\tiny \left(\begin{array}{c} 2 \ 2 \ 2 \ 4 \ 3 \\ 3 \ 3 \ 3 \ 3 \ 1 \end{array}\right)} & \cdot & (vii) \ (y,w)\in\DE \ \mbox{and} \ (y,x)\notin\DE \\ \cline{2-5}
 & x \leftarrow y & {\tiny \left(\begin{array}{c} 2 \ 2 \ 2 \ 3 \ 4 \\ 3 \ 3 \ 3 \ 4 \ 0 \end{array}\right)} & \checkmark & \Mp=\{w,x,y\},\ \Mm=\{w,x\} \\ \cline{2-5}
 & x \leftrightarrow y & {\tiny \left(\begin{array}{c} 2 \ 2 \ 2 \ 4 \ 4 \\ 3 \ 3 \ 3 \ 4 \ 1 \end{array}\right)} & \checkmark & \Mp=\{w,x,y\},\ \Mm=\{w,x\} \\ \hline\hline

\multirow{4}{*}{$[\C,\C]$} & x \cdots y & {\tiny \left(\begin{array}{c} 3 \ 3 \ 3 \ 3 \ 3 \\ 3 \ 3 \ 3 \ 3 \ 3 \end{array}\right)} & \cdot & (vii) \ (y,w)\in\DE \ \mbox{and} \ (y,x)\notin\DE \\ \cline{2-5}
 & x \rightarrow y & {\tiny \left(\begin{array}{c} 3 \ 3 \ 3 \ 4 \ 3 \\ 3 \ 3 \ 3 \ 3 \ 4 \end{array}\right)} & \cdot & (vii) \ (y,w)\in\DE \ \mbox{and} \ (y,x)\notin\DE \\ \cline{2-5}
 & x \leftarrow y & {\tiny \left(\begin{array}{c} 3 \ 3 \ 3 \ 3 \ 4 \\ 3 \ 3 \ 3 \ 4 \ 3 \end{array}\right)} & \cdot & (vii) \ (w,y)\in\DE \ \mbox{and} \ (x,y)\notin\DE \\ \cline{2-5}
 & x \leftrightarrow y & {\tiny \left(\begin{array}{c} 3 \ 3 \ 3 \ 4 \ 4 \\ 3 \ 3 \ 3 \ 4 \ 4 \end{array}\right)} & \checkmark & \Mp=\Mm=\{w,x,y\} \\ \hline
\end{array}
\]
\begin{center}
\footnotesize
\scshape
Table~\ref{tab:R2}
\it (continued)
\end{center}

\end{table}


\bibliography{hhchoice,hhchoice_books}

\begin{thebibliography}{10}

\bibitem{Berge:1973yq}
Berge.
\newblock {\em Graphs and Hypergraphs}.
\newblock North-Holland, Amsterdam, 1973.

\bibitem{Cameron:2007p7831}
Kathie Cameron, Elaine Eschen, Chinh Hoang, and R~Sritharan.
\newblock The complexity of the list partition problem for graphs.
\newblock {\em SIAM J. Discrete Math.}, 21(4):900--929, Jan 2007.

\bibitem{Erdos:1961p5414}
P~Erd{\H o}s and Tibor Gallai.
\newblock Graphs with prescribed degrees of vertices.
\newblock {\em Mat. Lapok}, 11:264--274, 1961.

\bibitem{Feder:2003p7761}
T~Feder, P~Hell, S~Klein, and R~Motwani.
\newblock List partitions.
\newblock {\em SIAM J. Discrete Math.}, Jan 2003.

\bibitem{Fulkerson:1960p3387}
Delbert~Ray Fulkerson.
\newblock Zero-one matrices with zero trace.
\newblock {\em Pacific J. Math.}, 10(3):831--836, Mar 1960.

\bibitem{hakimi:496}
S~Hakimi.
\newblock On realizability of a set of integers as degrees of the vertices of a
  linear graph. i.
\newblock {\em SIAM J. Appl. Math.}, 10(3):496--506, 1962.

\bibitem{hakimi:135}
S~Hakimi.
\newblock On realizability of a set of integers as degrees of the vertices of a
  linear graph ii. uniqueness.
\newblock {\em SIAM J. Appl. Math.}, 11(1):135--147, 1963.

\bibitem{Havel:1955p7653}
V~Havel.
\newblock A remark on the existence of finite graphs.
\newblock {\em {\v C}asopis Pest. Mat.}, 80:477--480, 1955.

\bibitem{Kleitman:1973p4546}
D~Kleitman and D~Wang.
\newblock Algorithms for constructing graphs and digraphs with given valences
  and factors.
\newblock {\em Discrete Math.}, 6(1):79--88, 1973.

\bibitem{McDonald:2007p6830}
John McDonald, Peter Smith, and Jonathan Forster.
\newblock Markov chain monte carlo exact inference for social networks.
\newblock {\em Social Networks}, 29(1):127--136, 2007.

\bibitem{Tripathi:2003p3520}
A~Tripathi and S~Vijay.
\newblock A note on a theorem of erd{\H o}s {\&} gallai.
\newblock {\em Discrete Math.}, 265(1-3):417--420, Apr 2003.

\end{thebibliography}
\bibliographystyle{plain}

\end{document}